\documentclass[11pt]{article}

\usepackage{tikz}
\usepackage{subfigure}
\usepackage{scalefnt}

\usepackage{amsmath,amsfonts,amsthm}

\usepackage{enumerate}

\newtheorem{defn}{{\bf Definition}}[section]
\newtheorem{eg}[defn]{{\bf Example}}
\newtheorem{lemma}[defn]{{\bf Lemma}}
\newtheorem{prop}[defn]{{\bf Proposition}}
\newtheorem{theo}[defn]{{\bf Theorem}}
\newtheorem{cor}[defn]{{\bf Corollary}}
\newtheorem{remark}[defn]{{\bf Remark}}
\newtheorem{conj}[defn]{{\bf Conjecture}}
\newtheorem{qn}[defn]{{\bf Question}}

\newcommand{\TPSS}{S^{\hspace{.2mm}2}\! \times \hspace{-3.3mm}_{-} \,
S^{\hspace{.1mm}1}}
\newcommand{\TPSSF}{S^{\hspace{.2mm}3}\! \times \hspace{-3.5mm}_{-} \,
S^{\hspace{.1mm}1}}
\newcommand{\TPSSD}{S^{\hspace{.2mm}d-1}\! \times \hspace{-3.5mm}_{-} \, S^{\hspace{.1mm}1}}
\newcommand{\TPSSDS}{S^{\hspace{.2mm}d-1}\! \times \hspace{-3.3mm}_{-} \, S^{\hspace{.1mm}1}}

\newcommand{\hoop}{-\hspace{-2mm}-\hspace{-1.5mm}\circ}
\newcommand{\hoops}{-\hspace{-2mm}-\hspace{-1.2mm}\circ}

\begin{document}

\title{\bf On stellated spheres and a tightness criterion for
combinatorial manifolds}
\author{{\bf Bhaskar Bagchi}$^{\rm a}$ and {\bf Basudeb
Datta}$^{\rm b}$ }

\date{}

\maketitle

\vspace{-5mm}

\noindent {\small $^{\rm a}$Theoretical Statistics and Mathematics Unit, Indian Statistical Institute,  Bangalore
560\,059, India.
\newline  \mbox{} \hspace{.2mm} bbagchi@isibang.ac.in

\smallskip

\noindent $^{\rm b}$Department of Mathematics, Indian Institute of Science, Bangalore 560\,012,  India. \newline
\mbox{} \hspace{.2mm}   dattab@math.iisc.ernet.in}


\begin{center}

\date{May 16, 2013}

\end{center}


\hrule

\bigskip

 \centerline{\sc Abstract}

\medskip

{\small We introduce the $k$-stellated spheres and consider the class ${\cal W}_k(d)$ of triangulated
$d$-manifolds all whose vertex links are $k$-stellated, and its subclass ${\cal W}^{\ast}_k(d)$ consisting of the
$(k+1)$-neighbourly members of ${\cal W}_k(d)$. We introduce the mu-vector of any simplicial complex and show
that, in the case of 2-neighbourly simplicial complexes, the mu-vector dominates the vector of its Betti numbers
componentwise; the two vectors are equal precisely for tight simplicial complexes. We are able to
estimate/compute certain alternating sums of the components of the mu-vector of any 2-neighbourly member of
${\cal W}_k(d)$ for $d\geq 2k$. As one consequence of this theory, we prove a lower bound theorem for such
triangulated manifolds, as well as determine the integral homology type of members of ${\cal W}^{\ast}_k(d)$ for
$d \geq 2k+2$. As another application, we prove that, when $d \neq 2k+1$, all members of ${\cal W}^{\ast}_k(d)$
are tight. We also characterize the tight members of ${\cal W}^{\ast}_k(2k + 1)$ in terms of their $k^{\rm th}$
Betti numbers. These results more or less answer a recent question of Effenberger, and also provide a uniform and
conceptual tightness proof for all except two of the known tight triangulated manifolds.

We also prove a lower bound theorem for triangulated manifolds in which the members of ${\cal W}_1(d)$ provide
the equality case. This generalises a result (the $d=4$ case) due to Walkup and K\"{u}hnel. As a consequence, it
is shown that every tight member of ${\cal W}_1(d)$ is strongly minimal, thus providing substantial evidence in
favour of a conjecture of K\"{u}hnel and Lutz asserting that tight triangulated manifolds should be strongly
minimal.}

\bigskip

{\small


\medskip

 }

\bigskip

\hrule

\section{Introduction}

All simplicial complexes considered here are finite and abstract. By a triangulated sphere/ ball/manifold, we
mean an abstract simplicial complex whose geometric carrier is a sphere/ ball/manifold. We identify two complexes
if they are isomorphic. Throughout, $\mathbb{F}$ is a fixed and arbitrary field. All homologies are simplicial
homologies with coefficients in $\mathbb{F}$; so we do not usually indicate the field in the notation for the
homology groups or Betti numbers.

A $d$-dimensional simplicial complex is called {\em pure} if all its maximal faces (called {\em facets}) are
$d$-dimensional. A $d$-dimensional pure simplicial complex is said to be a {\em weak pseudomanifold} if each of
its $(d - 1)$-faces is in at most two facets. For a $d$-dimensional weak pseudomanifold $X$, the {\em boundary}
$\partial X$ of $X$ is the pure subcomplex of $X$ whose facets are those $(d-1)$-dimensional faces of $X$ which
are contained in unique facets of $X$. The {\em dual graph} $\Lambda(X)$ of a pure simplicial complex $X$ is the
graph whose vertices are the facets of $X$, where two facets are adjacent in $\Lambda(X)$ if they intersect in a
face of codimension one. A {\em pseudomanifold} is a weak pseudomanifold with a connected dual graph. All
connected triangulated manifolds are automatically pseudomanifolds.

For any two simplicial complexes $X$ and $Y$, their {\em join} $X \ast Y$ is the simplicial complex whose faces
are the disjoint unions of the faces of $X$ with the faces of $Y$. (Here we adopt the convention that the empty
set is a face of every simplicial complex.)

For a finite set $\alpha$, let $\overline{\alpha}$ (respectively $\partial \alpha$) denote the simplicial complex
whose faces are all the subsets (respectively, all proper subsets) of $\alpha$. Thus, if $\#(\alpha)= n$,
$\overline{\alpha}$ is a copy of the standard triangulation $B^{\,n - 1}_n$ of the $(n- 1)$-dimensional ball, and
$\partial \alpha$ is a copy of the standard triangulation $S^{\,n-2}_n$ of the $(n-2)$-dimensional sphere. Thus,
for any two disjoint finite sets $\alpha$ and $\beta$, $\overline{\alpha}\ast\partial\beta$ and $\partial\alpha
\ast \overline{\beta}$ are two triangulations of a ball; they have identical boundaries, namely $(\partial
{\alpha}) \ast (\partial {\beta})$.

A subcomplex $Y$ of a simplicial complex $X$ is said to be an {\em induced} (or a {\em full}\,) subcomplex if
every face of $X$ contained in the vertex set of $Y$ is a face of $Y$. If $X$ is a $d$-dimensional simplicial
complex with an induced subcomplex $\overline{\alpha} \ast \partial \beta$ ($\alpha \neq \emptyset$, $\beta \neq
\emptyset$) of dimension $d$ (thus, $\dim(\alpha) + \dim(\beta) = d$), then $Y := (X\setminus (\overline{\alpha}
\ast \partial \beta)) \cup (\partial\alpha \ast \overline{\beta})$ is clearly another triangulation of the same
topological space $|X|$. In this case, $Y$ is said to be obtained from $X$ by the {\em bistellar move} $\alpha
\mapsto \beta$. If $\dim(\beta) = i$ ($0\leq i \leq d$), we say that $\alpha \mapsto \beta$ is a {\em bistellar
move of index $i$} (or an {\em $i$-move}, in short). Clearly, if $Y$ is obtained from $X$ by an $i$-move $\alpha
\mapsto \beta$ then $X$ is obtained from $Y$ by the (reverse) $(d- i)$-move $\beta \mapsto \alpha$. Notice that,
in case $i=0$, i.e., when $\beta$ is a single vertex, we have $\partial \beta = \{\emptyset\}$ and hence
$\overline{\alpha} \ast \partial \beta = \overline{\alpha}$. Therefore, our requirement that $\overline{\alpha}
\ast \partial \beta$ is the induced subcomplex of $X$ on $\alpha \sqcup \beta$ means that $\beta$ is a new
vertex, not in $X$. Thus, a $0$-move creates a new vertex, and correspondingly a $d$-move deletes an old vertex.
For $0 < i < d$, any $i$-move preserves the vertex set; these are sometimes called the {\em proper bistellar
moves}.

A triangulation $X$ of a manifold is called a {\em combinatorial manifold} if its geometric carrier $|X|$ is a
piecewise linear (pl) manifold with the pl structure induced from $X$. A combinatorial triangulation of a
sphere/ball  is called a {\em combinatorial sphere/ball} if it induces the standard pl structure (namely, that of
the standard sphere/ball) on its geometric carrier. Equivalently (cf. \cite{li, p}), a simplicial complex is a
combinatorial sphere (or ball) if it is obtained from a standard sphere (respectively, a standard ball) by a
finite sequence of bistellar moves. In general, a triangulated manifold is a combinatorial manifold if and only
if the link of each of its vertices is a combinatorial sphere or combinatorial ball. (Recall that the {\em link}
of a vertex $x$ in a complex $X$, denoted by ${\rm lk}_X(x)$, is the subcomplex $\{\alpha \in X \, : \, x\not \in
\alpha, \alpha \sqcup\{x\} \in X\}$.) This leads us to introduce\,:

\begin{defn} \label{stellated-sphere}
{\rm For $0\leq k \leq d+1$, a $d$-dimensional simplicial complex $X$ is said to be {\em $k$-stellated} if $X$
may be obtained from $S^{\,d}_{d+2}$ by a finite sequence of bistellar moves, each of index $< k$. By convention,
$S^{\,d}_{d+2}$ is the only $0$-stellated simplicial complex of dimension $d$.

Clearly, for $0\leq k \leq l \leq d+1$, $k$-stellated implies $l$-stellated. All $k$-stellated simplicial
complexes are combinatorial spheres. By Pachner's theorem \cite{p}, the $(d+ 1)$-stellated $d$-spheres are
precisely the combinatorial $d$-spheres.
 }
\end{defn}

We also recall\,:

\begin{defn} \label{k-stacked-ball}
{\rm For $0\leq k \leq d+1$, a triangulated $(d+1)$-dimensional ball $B$ is said to be {\em $k$-stacked} if all
the faces of $B$ of codimension (at least) $k+1$ lie in its boundary; i.e., if ${\rm skel}_{d - k}(B) = {\rm
skel}_{d - k}(\partial B)$. A triangulated $d$-sphere $S$ is said to be {\em $k$-stacked} if there is a
$k$-stacked $(d+ 1)$-ball $B$ such that $\partial B = S$.  }
\end{defn}

By an induction on the minimum number of bistellar moves needed (to obtain a $k$-stellated sphere beginning with
the standard sphere), it is easy to see that (cf. Propositions 2.7 and 2.9 in \cite{bd17})\,:

\begin{prop} \label{P1.3}
For $d\geq 2k -1$, a triangulated $d$-sphere $S$ is $k$-stellated if and only if $S$ is the boundary of a
shellable $k$-stacked $(d+1)$-ball. In consequence, all $k$-stellated spheres of dimension $\geq 2k-1$ are
$k$-stacked.
\end{prop}

As a result, any of the known examples of $k$-stacked non-shellable balls of dimension $\geq 2k$ lead to examples
of $k$-stacked spheres which are not $k$-stellated. For explicit examples, see Examples 3.1 in \cite{bd17}. Now
we introduce\,:

\begin{defn} \label{D1.4}
{\rm For $0 \leq k \leq d$, ${\cal W}_k(d)$ consists of the connected simplicial complexes of dimension $d$ all
whose vertex links are $k$-stellated $(d-1)$-spheres, and ${\cal K}_k(d)$ consists of the connected simplicial
complexes of dimension $d$ all whose vertex links are $k$-stacked $(d-1)$-spheres. }
\end{defn}

Thus, members of ${\cal W}_k(d)$ are combinatorial manifolds. As an immediate consequence of Proposition
\ref{P1.3}, we have\,:

\begin{cor} \label{C1.5}
${\cal W}_k(d) \subseteq {\cal K}_k(d)$ for $d \geq 2k$.
\end{cor}

It is easy to see that ${\cal W}_1(d) = {\cal K}_1(d)$. The class ${\cal K}_1(d)$ is called the {\em Walkup
class} (cf. \cite{wa}). Accordingly, the classes ${\cal K}_k(d)$ are known as the {\em generalised Walkup
classes} (cf. \cite{ef}).

The entire $g$-vector (equivalently, face vector) of any member of ${\cal K}_k(d)$ is determined by the $(k+1)$
numbers $g_1, \dots, g_{k+1}$. In particular, when $d\geq 2k$ is even, the Euler characteristic $\chi$ of any
member of ${\cal K}_k(d)$ is determined by the $(k+1)^{\rm st}$ component $g_{k+1}$ of its $g$-vector by the
formula (cf. \cite[Proposition 4.8]{bd13v2})
$$
(-1)^k{d+2 \choose k+1}(\chi - 2) =  2\,g_{k+1}.
$$

Recall that a simplicial complex $X$ is said to be {\em $l$-neighbourly} if any $l$ vertices of $X$ form a face
of $X$. If $M \in {\cal W}_k(d)$ is 2-neighbourly, with Betti numbers $\beta_i$ and $g$-components $g_i$, then we
show that
\begin{enumerate}
\item[{\rm (a)}] $g_{l+1} \geq {\displaystyle {d+2 \choose l+1} \sum_{i=1}^{l}} (-1)^{l-i}\beta_i$ \, for \, $1
\leq l < k$, \, when $d\geq 2k$, \vspace{-2mm}
\item[{\rm (b)}] $g_{k+1} \geq {\displaystyle {d+2 \choose k+1} \sum_{i=1}^{k}} (-1)^{k-i}\beta_i$, \, when $d =
2k+1$,
\item[{\rm (c)}] $g_{k+1} = {\displaystyle {d+2 \choose k+1} \sum_{i=1}^{k}} (-1)^{k-i}\beta_i$, \, when $d\geq
2k+2$, and
\item[{\rm (d)}] $\beta_i = 0$ \, for \, $k+1 \leq i \leq d-k-1$, \, when $d \geq 2k+2$.
\end{enumerate}
Since the components of the face vector are non-negative linear combinations of the $g$-numbers, this result may
be interpreted as a lower bound theorem for 2-neighbourly members of ${\cal W}_k(d)$.

It may be noted that Novik and Swartz have proved (\cite[Inequality (9)]{ns2}) that the inequalities (a), (b)
above hold for any triangulated $d$-manifold all whose vertex links are homology spheres satisfying the HLP (hard
Lefschetz property). Since a strong version of the GLBC (generalised lower bound conjecture) affirms that all
homology spheres have the HLP, it is expected that the inequalities (a) and (b) hold for all triangulated closed
$d$-manifolds (without the assumption of 2-neighbourliness or of $k$-stellated links). We conjecture that parts
(c) and (d) hold for any member of ${\cal K}_k(d)$ (without the 2-neighbourliness assumption). Indeed, these
considerations lead us to\,:

\begin{conj}[GLBC for triangulated manifolds] \label{Conj1.6}
Let $M$ be a connected triangulated closed $d$-manifold. Then, for 
$1\leq l \leq \frac{d- 1}{2}$, the $g$-numbers of $M$ satisfy 
$g_{l+1}(M) \geq \binom{d+2}{l+1} \sum_{i=1}^l (-1)^{l-i} \, 
\beta_i(M)$. Further, equality holds here for some $l
< \frac{d-1}{2}$ if and only if $M \in {\cal K}_l(d)$.
\end{conj}

Notice that Theorem \ref{P21} below proves the inequality of Conjecture \ref{Conj1.6} and the `if' part of the
equality case of the conjecture under the extra assumption that $M$ is 2-neighbourly and all the vertex links of
$M$ are $\lfloor \frac{d -1}{2}\rfloor$-stellated. The `$l=1$' case of Conjecture \ref{Conj1.6} (with $\mathbb{F}
= \mathbb{Q}$) was a conjecture of Kalai \cite{ka}; Novik and Swartz proved it (for any field $\mathbb{F}$, with
the extra hypothesis that $M$ is $\mathbb{F}$-orientable) in \cite[Theorem 5.2]{ns}. Since Conjecture
\ref{Conj1.6} includes the GLBC for homology spheres, we do not expect it to be settled in a hurry. The
inequalities in Conjecture \ref{Conj1.6} would follow from the strong GLBC quoted above.

We also prove a lower bound theorem (Theorem \ref{P23}) for general triangulated closed manifolds. When $d\geq
4$, among all triangulated closed $d$-manifolds with given first Betti number and given number of vertices,
members of ${\cal W}_1(d)$ (when they exist) minimize the face vector componentwise. The $d=4$ case of this
result is due to Walkup and K\"{u}hnel.

Using the Dehn-Sommerville equations, it is easy to see that a $k$-stellated/$k$-stacked sphere $S$ can be at
most $k$-neighbourly, unless $S$ is a standard sphere (cf. \cite[Corollary 4.6]{bd13v2}). Therefore a member of
${\cal W}_k(d)$/${\cal K}_k(d)$ can be at most $(k+ 1)$-neighbourly, unless it is a standard sphere. Let ${\cal
W}^{\ast}_k(d)$ (respectively ${\cal K}^{\ast}_k(d)$) be the class of $(k+ 1)$-neighbourly members of  ${\cal
W}_k(d)$ (respectively of ${\cal K}_k(d)$). The Dehn-Sommerville equations also imply that a $k$-neighbourly
$d$-sphere must be a standard $d$-sphere unless $d\geq 2k-1$. It follows that the classes ${\cal W}^{\ast}_k(d)$
and ${\cal K}^{\ast}_k(d)$ have no members except the standard $d$-sphere unless $d\geq 2k$. These remarks
indicate the importance of the classes ${\cal W}^{\ast}_k(d)$ and ${\cal K}^{\ast}_k(d)$. In conjunction with
Proposition \ref{P1.3} above, they also show that

\begin{cor} \label{C1.8}
${\cal W}_k^{\ast}(d) \subseteq {\cal K}_k^{\ast}(d)$ for all $d$ and $k$.
\end{cor}

We expect that part (b) of Theorem \ref{P23} should generalise as follows (compare K\"{u}hnel's conjecture
\cite[Conjecture 18]{lu2})\,:

\begin{conj} \label{Conj1.7}
If $M$ is an $m$-vertex connected triangulated closed $d$-manifold with Betti numbers $\beta_i$ then
${\,m+l-d-2\, \choose l+1} \geq {\,d+2\, \choose l+1} { \sum_{i=1}^l}(-1)^{l-i}\beta_i$ for $1 \leq l \leq
(d-1)/2$. Also, equality holds here for some $\,l< (d-1)/2$ if and only if $M \in {\cal K}_l^{\ast}(d)$.
\end{conj}

Moreover, we show that when $d \geq 2k+2$ and $k \geq 2$, any member of ${\cal W}^{\ast}_k(d)$ has the same
integral homology as the connected sum of $\beta$ copies of $S^{\,k} \times S^{\,d-k}$, where the non-negative
integer $\beta$ is given by the formula
$$
{{m+k-d-2 \choose k+1}} = {{d+2 \choose k+1}}\beta,
$$
with $m = f_0$, the number of vertices. This result may be compared with Kalai's theorem\,: for $d \geq 4$, any
member of ${\cal W}_1(d)$ triangulates the connected sum of finitely many copies of $S^1 \times S^{\,d-1}$ or
$\TPSSD$.

Recall that a connected simplicial complex $X$ is said to be {\em $\mathbb{F}$-tight} if the inclusion map from
any induced subcomplex of $X$ into $X$ is injective at the level of $\mathbb{F}$-homology. In case of a
triangulated closed manifold $X$, this has the following geometric interpretation\,: $X$ is $\mathbb{F}$-tight if
the standard geometric realization of $X$ in $\mathbb{R}^{\,n-1}$ ($n=$ number of vertices of $X$) is ``as convex
as possible" subject to the constraint imposed by its homology with $\mathbb{F}$-coefficients. Since any induced
subcomplex of a tight simplicial complex is obviously tight, tightness imposes an extremely powerful constraint
on the possible combinatorics of a simplicial complex. For instance, any $\mathbb{F}$-tight simplicial complex is
necessarily 2-neighbourly and any $\mathbb{F}$-tight triangulated closed manifold is $\mathbb{F}$-orientable.
Thus, it is not surprising that, apart from four infinite families (including the trivial family of standard
spheres), only thirty sporadic examples of tight triangulated manifolds (of dimension $>2$) are known so far.

All this makes it very important to obtain usable combinatorial criteria for tightness of triangulated manifolds.
In this paper, we introduce the sigma- and mu-vector (with respect to $\mathbb{F}$) of simplicial complexes. The
sigma-vector of a simplicial complex $X$ is an weighted average of the beta-vectors (i.e., the vectors of Betti
numbers) of the induced subcomplexes of $X$. The mu-vector of $X$ is essentially the average of the sigma-vectors
of the vertex links of $X$. In \cite{k95}, K\"uhnel uses the notion of regular simplexwise linear (rsl) real
valued functions on the vertex set of a simplicial complex $X$. With any such function, one associates a vector
of Morse indices (with respect to $\mathbb{F}$) of the critical points of the function. These vectors are used in
\cite{k95} to investigate $\mathbb{F}$-tightness of $X$. One observes that, in this investigation, only the
linear order induced on the vertex set $V(X)$ of $X$ (as a pull back of the usual linear order on $\mathbb{R}$ by
the given rsl function) is important, and not the rsl functions themselves. So, we might say that two rsl
functions are essentially the same if they induce the same linear ordering on $V(X)$. Thus, if $X$ has $n$
vertices, then there are only $n!$ essentially distinct rsl functions on $X$. Lemma \ref{L5} below shows that,
when $X$ is 2-neighbourly, the mu-vector of $X$ as defined here is just the average of the vectors of
Morse-indices over the $n!$ essentially distinct rsl functions. So, when $X$ is 2-neighbourly, K\"uhnel's
criteria for tightness of $X$ may be reformulated in terms of its mu-vector alone. This is done in Theorem
\ref{P18} below. This theorem is an immediate consequence of Theorem \ref{P16} which is our version of K\"uhnel's
combinatorial Morse theory. It shows that for a 2-neighbourly simplicial complex $X$, the investigation of
$\mathbb{F}$-tightness of $X$ boils down to an investigation of the sigma-vectors of the vertex links of $X$. We
also observe (Proposition \ref{P17} below) that, to be $\mathbb{F}$-tight, $X$ must be 2-neighbourly. Thus, the
restriction to 2-neighbourly simplicial complexes in Theorem \ref{P16} is no real loss of generality. In Theorem
\ref{P19} below, we estimate/compute certain alternating sums of the sigma-components of $k$-stellated spheres of
dimension $\geq 2k-1$. Apart from the introduction of the sigma- and mu-vectors, this theorem may be regarded as
the main contribution of this paper. As an immediate consequence of Theorem \ref{P19}, we obtain (Theorem
\ref{P20}) estimation/computation of the corresponding alternating sums of the mu-components of 2-neighbourly
members of ${\cal W}_k(d)$ for $d\geq 2k$. In conjunction with the strong Morse inequalities of Theorem \ref{P16}
(a), Theorem \ref{P20} implies the lower bound theorem for 2-neighbourly members of ${\cal W}_k(d)$ (Theorem
\ref{P21}) as well as the determination (Theorem \ref{P24}) of the integral homotopy type of members of ${\cal
W}^{\ast}_k(d)$, $d\geq 2k+2$.

In \cite{ef}, Effenberger asked whether, for $d\neq 2k+1$, all members of ${\cal K}^{\ast}_k(d)$ are
$\mathbb{F}$-tight. As another application of Theorem \ref{P20}, we give a partial answer to this question by
showing (Theorem \ref{P25}) that, for $d\neq 2k+1$, all $\mathbb{F}$-orientable members of  ${\cal
W}^{\ast}_k(d)$ are $\mathbb{F}$-tight. (To be $\mathbb{F}$-tight, a triangulated closed manifold must be
$\mathbb{F}$-orientable, cf. Proposition \ref{P17} below. Also note that, when $k\geq 2$, members of  ${\cal
K}^{\ast}_k(d)$ are necessarily simply connected and hence $\mathbb{F}$-orientable.) In \cite{bd15}, we show that
all polytopal upper bound spheres (for instance cyclic spheres) of dimension $2k+1$ belong to the class ${\cal
W}^{\ast}_k(2k+1) \subseteq {\cal K}^{\ast}_k(2k+1)$. Since it is easy to see (cf. Remark \ref{R1} below) that
the standard $d$-spheres are the only $\mathbb{F}$-tight triangulated spheres, the examples of polytopal upper
bound spheres show that the hypothesis $d\neq 2k+1$ in Theorem \ref{P25} is essential. We also obtain a necessary
and sufficient condition for members of ${\cal W}^{\ast}_k(2k+1)$ to be $\mathbb{F}$-tight, involving its $k^{\rm
th}$ Betti number. Using these results, we are able to present a conceptual tightness proof of all except two of
the known tight triangulated closed manifolds.

We recall that a simplicial complex $X$ is said to be {\em minimal} if it has the fewest number of vertices among
all possible triangulations of the geometric carrier $|X|$ of $X$. Let's say that $X$ is {\em strongly minimal}
if it achieves the componentwise smallest face vector among all triangulations of $|X|$. Our interest in
tightness stems from the following conjecture of K\"{u}hnel and Lutz \cite{kl}.

\begin{conj}[K\"{u}hnel and Lutz] \label{Conj1.9}
Every $\mathbb{F}$-tight simplicial complex is strongly minimal.
\end{conj}

As a consequence of Theorem \ref{P23}, we have (Corollary \ref{C11} below) that all the tight members of ${\cal
W}^{\ast}_1(d) = {\cal K}^{\ast}_1(d)$ are strongly minimal. This result may be regarded as a substantial
evidence in favour of Conjecture \ref{Conj1.9}.

\section{The mu-vector, Morse theory and tightness}

\noindent {\bf Notation\,: (a)} For any simplicial complex $X$ with vertex set $V(X)$ and any set $A$, $X[A]$
will denote the induced subcomplex of $X$ on the vertex set $A\cap V(X)$. Thus, $X[A] = \{\alpha\in X \, :\,
\alpha \subseteq A\}$. \newline {\bf (b)} For any set $V$ and integer $i\geq 0$, we denote by ${V \choose i}$ the
collection of all subsets of $V$ of size $i$. Thus, ${V \choose i}= \{A\subseteq V \, : \, \#(A)=i\}$.

\bigskip

Recall that all homologies are with respect to the field $\mathbb{F}$. Now, we introduce

\begin{defn} \label{beta-vector}
{\rm  Let $X = X^d_m$ be a simplicial complex of dimension $d$ on $m$ vertices. Then we define the beta-, sigma-,
and mu-vector of $X$ as follows. The {\em beta-vector} of $X$ is the vector $(\beta_0, \beta_1, \dots, \beta_d)$,
where $\beta_i=\beta_i(X)$ is the $i^{\rm th}$ Betti number of $X$. That is, $\beta_i(X) =
\dim_{\mathbb{F}}H_i(X)$, $0\leq i\leq d$. More generally, if $Y$ is a subcomplex of $X$, we use $\beta_i(X, Y)$
to denote $\dim_{\mathbb{F}}H_i(X, Y)$, $0\leq i\leq d$. As usual, $\widetilde{\beta}_i$ will denote the
corresponding reduced Betti numbers. Thus, $\widetilde{\beta}_i = \beta_i$ if $i \neq 0$ and
 $\widetilde{\beta}_0 = \beta_0-1$.

The {\em sigma-vector} $(\sigma_0, \sigma_1, \dots, \sigma_d)$ of $X$ is defined by
$$
\sigma_i = \sigma_i(X) = \sum_{j=0}^m \frac{1}{{m\choose j}} \sum_{A\in {V(X) \choose j}}
\widetilde{\beta}_i(X[A]), ~~ 0\leq i\leq d.
$$

We define the {\em mu-vector} $(\mu_0, \dots, \mu_d)$ of $X$ by
\begin{eqnarray*}
\mu_0 & = & \mu_0(X) = 1 \\
\mu_i & = & \mu_i(X) = \delta_{i1} + \frac{1}{m}\sum_{x \in V(X)} \sigma_{i-1}({\rm lk}_X(x)), ~~ 1 \leq i\leq d.
\end{eqnarray*}
(By a slight abuse of notation, let's write $\emptyset$ for the trivial simplicial complex whose only face is the
empty set. Here, we have adopted the convention $\widetilde{\beta}_0(\emptyset) = - 1$ and
$\widetilde{\beta}_i(\emptyset) = 0$ if $i \neq 0$. This convention accounts for the Kronecker delta in the
definition of the mu-vector.) }
\end{defn}

\begin{lemma} \label{L2.2}
Let $S$ be a triangulated homology $d$-sphere. Then the sigma-vector of $S$ satisfies
$$
\sigma_{d-1-i} = \left\{
\begin{array}{ll}
\sigma_i & \mbox{ if \, $0<i<d-1$} \\
\sigma_i +1 & \mbox{ if \, $i = 0$} \\
\sigma_i -1  &  \mbox{ if \, $i = d-1$}.
\end{array}
\right.
$$
\end{lemma}

\noindent {\bf Proof.} Let $V$ be the vertex set of $S$. For any subset $A$ of $V$, let $\bar{A} = V\setminus A$.
Then, as an immediate consequence of Alexander duality and the exact sequence for pairs, we have
$\widetilde{\beta}_{d-1-i}(S[A]) = \widetilde{\beta}_{i}(S[\bar{A}])$ for $0\leq i\leq d-1$, except when $A=V$,
$i=0$ or $A=\emptyset$, $i=d-1$. Taking the appropriate weighted average of these equations over all $A$ (and
noting that $A \mapsto \bar{A}$ is a bijection on the power set of $V$), we get the result. \hfill $\Box$

\begin{theo} \label{T2.3}
The mu-vector of any triangulated closed $d$-manifold satisfies $\mu_{d-i} = \mu_{i}$, $0\leq i\leq d$.
\end{theo}

\noindent {\bf Proof.} Let $M$ be a triangulated closed $d$-manifold. For any vertex $x$ of $M$, we apply Lemma
\ref{L2.2} to the link of $x$ (with $d-1$ in place of $d$ and $i-1$ in place of $i$). When the result is averaged
over all vertices $x$ of $M$, we get $\mu_{d-i} = \mu_{i}$, $1\leq i\leq d-1$. Also, for any homology $d$-sphere
$S$, $\widetilde{\beta}_d$ of any proper induced subcomplex of $S$ is zero, while $\widetilde{\beta}_d$ of $S$ is
one. Hence we get $\sigma_d(S) = 1$. Applying this observation to the vertex links of $M$ and averaging, we get
$\mu_d = 1 =\mu_0$. \hfill $\Box$

\bigskip

We are mildly surprised that (unlike the corresponding result for the vector of Betti numbers), this theorem does
not require the $\mathbb{F}$-orientability of the manifold.

\bigskip

\noindent {\bf Notation\,:} We write $\hoop$ for the covering relation for set inclusion. Thus, for sets $A$ and
$B$, $A \hoop B$ means that $A \subseteq B$ and $\#(B\setminus A) =1$.

\medskip

With this notation, we have\,:

\begin{lemma} \label{L5}
Let $X$ be a $2$-neighbourly simplicial complex of dimension $d$ on $m$ vertices. Then, the mu-vector of $X$
$($with respect to any field$)$ is given by\,: \vspace{-2mm}
$$
\mu_i= \frac{1}{m}\sum_{j=1}^m \frac{1}{{m-1\choose j-1}} \sum_{\scriptsize
\begin{array}{c}
A, B \subseteq V(X),\\
\#(B) = j, \\ A \hoops B
\end{array}
} \hspace{-4mm}\beta_i(X[B], X[A]), ~~~ 0\leq i \leq d.
$$
\end{lemma}

\noindent {\bf Proof.} Since $X$ is 2-neighbourly, so is any induced subcomplex of $X$. Hence, for $A \hoop B
\subseteq V(X)$, $\beta_0(X[B], X[A]) =1$ if $\#(B)=1$, $A = \emptyset$ and $\beta_0(X[B], X[A]) =0$ otherwise.
Hence the formula holds for $i=0$. So, let $1\leq i \leq d$.  For any $x\in V(X)$, let $L_x$ be the link of $x$
in $X$. Then, each $L_x$ is a simplicial complex of dimension $d-1$ with exactly $m-1$ vertices. Let $V= V(X)$ be
the vertex set of $X$. Thus, the vertex set of $L_x$ is $V \setminus\{x\}$. We have
$$
\sigma_{i-1}(L_x) = \sum_{j=1}^m \frac{1}{{m-1 \choose j-1}} \sum_{A \in {V\setminus\{x\} \choose j-1}}
\widetilde{\beta}_{i -1}(L_x[A]).
$$
But the exact homology sequence for pairs and the excision theorem yield
$$
\widetilde{\beta}_{i-1}(L_x[A]) = \left\{
\begin{array}{ll}
\beta_i(x \ast L_x[A], L_x[A]) = \beta_i(X[A \sqcup\{x\}], X[A]) &
\mbox{ when } ~  A \neq \emptyset ~ \mbox{ and}\\
 \beta_i(X[A \sqcup\{x\}], X[A])
- \delta_{i1} & \mbox{ when } ~  A = \emptyset.
\end{array}
\right.
$$
Therefore, from the definition of the mu-vector of $X$, we have
$$
\mu_i= \frac{1}{m}\sum_{x\in V}\sum_{j=1}^m \frac{1}{{m-1\choose j-1}} \sum_{A \in {V\setminus\{x\} \choose j-1}}
\beta_i(X[A\sqcup\{x\}], X[A]).
$$ \hfill $\Box$

\bigskip

The following linear algebra lemma must be well known. But, we could not find a reference to it in the required
form.

\begin{lemma} \label{L6}
Let $V_1 \stackrel{T_1}{\longrightarrow} V_2 \rightarrow \cdots \rightarrow V_{2m}
\stackrel{T_{2m}}{\longrightarrow}  V_{2m+1}$ be an exact sequence of linear transformations between finite
dimensional vector spaces $($involving an even number $2m$ of arrows$)$. Then $\sum_{i=1}^{2m+1} (-1)^i \dim(V_i)
\leq 0$. Equality holds here if and only if $T_1$ is injective and $T_{2m}$ is surjective.
\end{lemma}

\noindent {\bf Proof.} From the assumed exactness, we have $\dim(V_i) = {\rm rank}(T_i) + {\rm rank}(T_{i-1})$
for $1 < i < 2m+1$. Therefore, the negative of the alternating sum above telescopes to $\dim({\rm Ker}(T_1)) +
\dim({\rm Coker}(T_{2m}))$. Hence the result. \hfill $\Box$

\bigskip

The following theorem is our version of the usual combinatorial Morse theory. In particular, the parts (a) and
(b) of this theorem are the strong and weak Morse inequalities (averaged over all possible ``regular simplexwise
linear" functions; compare \cite{k95}). For an alternative combinatorial version of Morse theory, consult
\cite{fo}. As we shall see, the version developed here is specially suited to the study of $\mathbb{F}$-tightness
of $\mathbb{F}$-orientable 2-neighbourly triangulated closed manifolds.

\begin{theo}  \label{P16}
Let $X$ be a $2$-neighbourly simplicial complex of dimension $d$. Then the mu- and beta-vectors of $X$ are
related as follows. \vspace{-1mm}
\begin{enumerate}
\item[{\rm (a)}] For $0\leq j \leq d$,\, $\sum_{i=0}^j (-1)^{j-i} \mu_i \geq \sum_{i=0}^j (-1)^{j-i} \beta_i$,
with equality for $j=d$.
 \vspace{-1mm}
\item[{\rm (b)}]  For $0\leq j \leq d$,\, $\mu_j \geq \beta_j$.
 \vspace{-1mm}
\item[{\rm (c)}] The following are equivalent for any fixed index $j$ $(0\leq j\leq d)$\,:
\begin{enumerate}
\item[$(i)$] $\sum_{i=0}^j (-1)^{j-i} \mu_i = \sum_{i=0}^j (-1)^{j - i} \beta_i$, and \vspace{-1mm} \item[$(ii)$]
for any induced subcomplex $Y$ of $X$, the morphism $H_{j}(Y; \mathbb{F}) \to H_{j}(X; \mathbb{F})$ induced by
the inclusion map $Y \hookrightarrow X$ is injective.
\end{enumerate}
 \vspace{-1mm}
\item[{\rm (d)}] The following are equivalent for any fixed index $j$ $(0\leq j\leq d)$\,: \vspace{-1mm}
\begin{enumerate}
\item[$(iii)$] $\mu_j = \beta_j$, and
\item[$(iv)$] for any induced subcomplex $Y$ of $X$, both the morphisms $H_{j-1}(Y; \mathbb{F}) \to H_{j-1}(X;
\mathbb{F})$ and $H_{j}(Y; \mathbb{F}) \to H_{j}(X; \mathbb{F})$ induced by the inclusion map $Y \hookrightarrow
X$ are injective.
\end{enumerate}
 \vspace{-1mm}
\item[{\rm (e)}] If, further, $X$ is an $\mathbb{F}$-orientable closed manifold, then $\beta_{d-j} = \beta_{j}$
and $\mu_{d-j} = \mu_{j}$ for $0\leq j\leq d$.
\end{enumerate}
\end{theo}

\noindent {\bf Proof.} (a) Fix an index $j$ and subsets $A \subseteq B$ of $V(X)$. We have the following exact
sequence of relative homology\,: $H_{j}(X[A]) \to H_{j}(X[B]) \to H_{j}(X[B], X[A]) \to H_{j-1}(X[A]) \to \cdots
\to H_{0}(X[A]) \to H_{0}(X[B]) \to H_{0}(X[B], X[A]) \to 0$. If necessary, we may append an extra $0\to 0$ at
the extreme right, to ensure that this exact sequence has an even number of arrows. Applying Lemma \ref{L6} to
this sequence, we get\,:
\begin{eqnarray} \label{eq8}
\sum_{i=0}^j (-1)^{j-i}\beta_i(X[B], X[A]) & \geq & \sum_{i=0}^j (-1)^{j-i}\left(\beta_i(X[B]) -
\beta_i(X[A])\right)
\end{eqnarray}
for all pairs $A \subseteq B$ of subsets of $V(X)$.

Since the extreme right arrow in the above sequence is trivially a surjection, Lemma \ref{L6} says that, for any
given pair $A\subseteq B$, equality in (\ref{eq8}) holds if and only if the morphism $H_{j}(X[A]) \to
H_{j}(X[B])$ induced by the inclusion map $A \hookrightarrow B$ is an injection.

Now, in view of Lemma \ref{L5}, taking the appropriate weighted sum of the inequalities (\ref{eq8}) over all
pairs $(A, B)$ with $A \hoop B$, we get
\begin{eqnarray} \label{eq9}
\sum_{i=0}^j (-1)^{j-i}\mu_i & \geq & \frac{1}{m} \sum_{l=1}^m \frac{1}{{m-1 \choose l-1}} \sum_{\scriptsize
\begin{array}{c}
A, B \subseteq V(X),\\
\#(B) = l, \\ A \hoops B
\end{array}
} \sum_{i=0}^j(-1)^{j-i}\left(\beta_i(X[B]) - \beta_i(X[A])\right).
\end{eqnarray}
Here $m = \#(V(X))$.

Equality holds in (\ref{eq9}) if and only if $H_j(X[A]) \to H_j(X[B])$ is injective for all pairs $(A, B)$ with
$A \hoops B\subseteq V(X)$. In particular, since $X$ is $d$-dimensional and each $d$-cycle of $X[A]$ is a
$d$-cycle of $X[B]$, equality holds in (\ref{eq9}) for $j = d$. The right hand side of (\ref{eq9}) may be written
as $\sum_{i= 0}^j (-1)^{j-i} \sum_{C \subseteq V(X)} \alpha(C) \beta_i(X[C])$, where the coefficients $\alpha(C)$
are given in terms of $n :=\#(C)$ by the formula
$$
\alpha(C) \, = \, \frac{1}{m} \frac{1}{{m-1 \choose n-1}} \sum_{\scriptsize
\begin{array}{c}
A \subseteq V(X),\\
 A \hoops C
\end{array}}\hspace{-4mm}1 ~
- ~ \frac{1}{m} \frac{1}{{m-1 \choose n}} \sum_{\scriptsize
\begin{array}{c}
B \subseteq V(X),\\
C \hoops B
\end{array}}\hspace{-4mm}1
 \, = \, \frac{1}{m} \left(\frac{n}{{m-1 \choose n-1}} -
\frac{m-n}{{m-1 \choose n}}\right),
$$
where the first term occurs only for $n>0$, and the second term occurs only for $n<m$. This simplifies to
$$
\alpha(C) = \left\{\begin{array}{rl}
+1 & \mbox{if } ~ C =V(X) \\
-1 & \mbox{if } ~ C =\emptyset \\
0 & \mbox{otherwise}.
\end{array}
\right.
$$
Therefore, the right hand side of (\ref{eq9}) simplifies to ${\displaystyle \sum_{i=0}^j(-1)^{j-i}(\beta_i(X) -
\beta_i(\emptyset)) = \sum_{i=0}^j(-1)^{j-i}\beta_i.}$

(b) We have,
$$
\mu_j = \sum_{i=0}^{j-1} (-1)^{j-1-i} \mu_i + \sum_{i=0}^{j} (-1)^{j-i} \mu_i \geq \sum_{i=0}^{j-1} (-1)^{j-1-i}
\beta_i + \sum_{i=0}^{j} (-1)^{j-i} \beta_i = \beta_j.
$$

(c) By the proof of part (a), we see that the equality $(i)$ holds if and only if $H_{j}(X[A]) \to H_{j}(X[B])$
is injective for all pairs $A \hoop B \subseteq V(X)$. Now, let $Y$ be an induced subcomplex of $X$, say with
vertex set $A$. Take a sequence $A = A_1 \hoop A_2 \hoop \cdots \hoop A_n = V(X)$. Then the morphism $H_{j}(Y)
\to H_{j}(X)$ is the composition of the morphisms $H_{j}(X[A_1]) \to H_{j}(X[A_2]) \to \cdots \to H_{j}(X[A_n])$.
If $(i)$ holds then $H_{j}(Y) \to H_{j}(X)$, being a composition of injective morphisms, is itself injective.
Thus, $(ii)$ holds. Conversely, if $(ii)$ holds, then for any pair $A_1 \hoop A_2 \subseteq V(X)$, choose
$Y=X[A_1]$. Then the composition of the above sequence of morphisms is injective. So, the first morphism
$H_{j}(X[A_1]) \to H_{j}(X[A_2])$ in this sequence must be injective. Therefore $(i)$ holds.

(d) From the proof of part (b), we see that $(iii)$ holds if and only if $\sum_{i=0}^{j-1} (-1)^{j-1-i} \mu_i =
\sum_{i=0}^{j-1} (-1)^{j-1-i} \beta_i$ and $\sum_{i=0}^{j} (-1)^{j-i} \mu_i = \sum_{i=0}^{j} (-1)^{j-i} \beta_i$.
Therefore, part (d) follows from part (c).

(e) In this case, Poincar\'e duality yields $\beta_{d-i}= \beta_i$. Theorem \ref{T2.3} yields $\mu_{d-i} =
\mu_i$. \hfill $\Box$

\bigskip

Now we recall\,:

\begin{defn} \label{tight}
{\rm Let $X$ be a $d$-dimensional simplicial complex and $\mathbb{F}$ be a field. We say that $X$ is {\em tight
with respect to} $\mathbb{F}$ (or, in short, {\em $\mathbb{F}$-tight}) if (i) $X$ is connected, and (ii) for all
induced subcomplexes $Y$ of $X$ and for all $0\leq j \leq d$, the morphism $H_{j}(Y; \mathbb{F}) \to H_{j}(X;
\mathbb{F})$ induced by the inclusion map $Y \hookrightarrow X$ is injective.

Note that, for fields $\mathbb{F}_1 \subseteq \mathbb{F}_2$, $X$ is $\mathbb{F}_1$-tight if and only if $X$ is
$\mathbb{F}_2$-tight. Therefore, in studying $\mathbb{F}$-tightness, we may, without loss of generality, restrict
to prime fields $\mathbb{F}$, i.e., $\mathbb{F}=\mathbb{Q}$ or $\mathbb{F}=\mathbb{Z}_p$, $p$ prime. Moreover,
for any simplicial complex $X$, the following are equivalent\,: (a) $X$ is $\mathbb{F}$-tight for all fields
$\mathbb{F}$, (b) $X$ is $\mathbb{Z}_p$-tight for all primes $p$, and (c) $X$ is $\mathbb{Q}$-tight. In view of
this observation, we shall say that $X$ is {\em tight} if it is $\mathbb{Q}$-tight. }
\end{defn}

\begin{remark} \label{R1}
{\rm Clearly, if $X$ is $\mathbb{F}$-tight then so is every induced subcomplex of $X$. From this trivial
observation, it is easy to see that the standard sphere $S^{\hspace{.1mm}d}_{d + 2}$ is the only
$\mathbb{F}$-tight triangulated $d$-sphere, and the standard ball $B^{\hspace{.1mm}d}_{d +1}$ is the only
$\mathbb{F}$-tight triangulated $d$-ball. Thus, the standard spheres/balls are the only possible induced
spheres/balls in an $\mathbb{F}$-tight simplicial complex. Since the triangulated ball $\overline{\alpha}\ast
\partial\beta$ is a standard ball only if $\dim(\beta) = 0$, it follows that an $\mathbb{F}$-tight simplicial
complex does not allow any bistellar move of non-zero index. }
\end{remark}

The following result is well known (see \cite{kl} for example). We have included its short proof for
completeness, since we could not find the proof in the existing literature.

\begin{prop}  \label{P17}
Let $X$ be an $\mathbb{F}$-tight simplicial complex.
 \vspace{-1mm}
\begin{enumerate}
\item[{\rm (a)}] If $X$ is $(k-1)$-connected in the sense of homotopy $($for some $k \geq1)$ then $X$ is
$(k+1)$-neighbourly.
 \vspace{-1mm}
\item[{\rm (b)}] If $X$ is a triangulated closed manifold, then $X$ is $\mathbb{F}$-orientable.
\end{enumerate}
\end{prop}

\noindent {\bf Proof.} (a) Suppose not. Let $l$ be the smallest integer such that $X$ is not $(l+1)$-neighbourly.
We have $1\leq l \leq k$. The induced subcomplex $X[\alpha]$ of $X$ on any missing $l$-face $\alpha$ is an
$S^{\hspace{.2mm}l-1}_{l+1}$. Since $H_{l-1}(X[\alpha]; \mathbb{F}) \to H_{l-1}(X; \mathbb{F})$ is injective, it
follows that $\widetilde{H}_{l-1}(X; \mathbb{F}) \neq 0$. This is a contradiction since $l\leq k$ and $X$ is
$(k-1)$-connected.

(b) $\mathbb{F}$-orientability of $X$ follows since $\beta_d = \mu_d = \mu_0 = 1$, where the first equality is by
Theorem \ref{P16}\,(d), and the second equality is by Theorem \ref{T2.3}.  \hfill $\Box$

\bigskip

Another way of stating Proposition \ref{P17}\,(b) is that, if $X$ is a triangulated closed manifold which is not
orientable (over $\mathbb{Z}$), then $\mathbb{F}= \mathbb{Z}_2$ is essentially the only choice for a field for
which $X$ has a chance of being $\mathbb{F}$-tight. Note that, as a special case ($k=1$) of Proposition
\ref{P17}\,(a), any $\mathbb{F}$-tight simplicial complex is necessarily 2-neighbourly (whatever the field
$\mathbb{F}$).

Now we have\,:

\begin{theo}  \label{P18}
A simplicial complex $X$ of dimension $d$ is $\mathbb{F}$-tight if and only if $X$ is $2$-neighbourly and
$\mu_i(X; \mathbb{F}) = \beta_i(X; \mathbb{F})$ for all indices $i$, $0\leq i \leq d$.
\end{theo}

\noindent {\bf Proof.} This is immediate from Theorem \ref{P16}\,(d) and Proposition \ref{P17}\,(a). \hfill
$\Box$

\begin{remark} \label{R2}
{\rm One can see from the proof of Lemma \ref{L5} that the $i^{\rm th}$ component $\mu_i$ of the mu-vector of a
2-neighborly simplicial complex $X$ is the average (over all vertices of $X$ and all rsl functions on $X$) of the
usual Morse multiplicities in degree $i$. (But it would take too much effort and space to make this explicit.)
Accordingly, the Morse inequalities of Theorem \ref{P16} are not really new, but are averaged versions of the
well known combinatorial Morse inequalities of K\"{u}hnel \cite{k95}. Therefore, in principle, parts (a), (b) and
(e) of Theorem \ref{P16} could be deduced from the known Morse inequalities by averaging. However, in our
opinion, one important achievement of this paper is the clarification of exactly when equality holds in (the
averaged version of) these inequalities, namely, parts (c) and (d) of Theorem \ref{P16}. In order to establish
parts (c) and (d), we found it necessary to reproduce the (essentially known) proof of parts (a) and (b) as well
-- in the language of the mu-vector. The introduction of the single mu-vector (in place of the entire families of
rsl functions) also clarifies the relationship (Theorem \ref{P18}) between these notions and tightness. For these
reasons, as well as for the benefit of new researchers in the field, we decided to include the self-contained
proof of Theorem \ref{P16}, even though some parts of this proof may sound familiar to the experts.}
\end{remark}

\section{The main results}

For a $d$-dimensional simplicial complex $X$, $f_i = f_i(X)$ denotes the number of $i$-dimensional faces of $X$
($-1 \leq i\leq d$). Thus, $f_{-1} = 1$, corresponding to the empty face of $X$. The vector {\boldmath $f$}$(X)
=(f_0, f_1, \dots, f_d)$ is called the {\em face vector} (or {\em $f$-vector}) of $X$.

The $g$-vector $g(X)= (g_0, \dots, g_{d+1})$ of a $d$-dimensional simplicial complex $X$ is defined in terms of
its $f$-vector by\,:
\begin{equation} \label{eq1}
g_j(X) = \sum_{i=-1}^{j-1} (-1)^{j-i-1}{d-i+1 \choose j-i-1}f_i(X), ~~ 0 \leq j \leq d+1.
\end{equation}
(\ref{eq1})  may be inverted to obtain the $f$-vector of $X$ in terms of its $g$-vector\,:
\begin{equation} \label{eq2}
f_i(X) = \sum_{j=0}^{i+1} {d-j+2 \choose i-j+1}g_j(X), ~~ -1 \leq i \leq d.
\end{equation}

The following interpretation of the $g$-vector is well known (cf. \cite[p. 83]{p}, for instance).

\begin{lemma} \label{L3}
Let $X$, $Y$ be two simplicial complexes of dimension $d$. If $Y$ is obtained from $X$ by a single bistellar move
of index $l$ then, for $0\leq j \leq d$,
$$
g_{j+1}(Y) - g_{j+1}(X) = \left\{
\begin{array}{rl}
+1 & \mbox{ if $j = l \neq d/2$} \\
-1 & \mbox{ if $j = d-l \neq d/2$} \\
0  & \mbox{ otherwise}.
\end{array}
\right.
$$
\end{lemma}

The following lemma is also well known (see \cite[Proposition 2.3]{s04}, for example).

\begin{lemma} \label{L4}
If $X$ is a simplicial complex of dimension $d$, then
$$
\sum_{x \in V(X)} g_j({\rm lk}_X(x)) = (d+2-j)g_j(X) + (j+1) g_{j+1}(X) ~ \mbox{ for } ~ 0\leq j \leq d.
$$
\end{lemma}

In the following result, $B^{\hspace{.2mm}e}$ stands for an arbitrary triangulated ball of dimension $e$, and
$S^{\hspace{.2mm}e -1}$ stands for an arbitrary triangulated sphere of dimension $e-1$. Here $e \geq 0$, and
$S^{-1} = \emptyset$.

\begin{lemma} \label{L7}
Let $X$, $Y$ be simplicial complexes such that $Y = X \cup B^{\hspace{.2mm}e}$. Suppose $X \cap
B^{\hspace{.2mm}e} = B^{\hspace{.2mm}e-1}$ $($with $e\geq 1)$ or $X \cap B^{\hspace{.1mm}e} =
S^{\hspace{.2mm}e-1}$ $($with $e\geq 0)$. Then the reduced Betti numbers $($with respect to any field
$\mathbb{F}\,)$ of $X$ and $Y$ are related as follows. \vspace{-1mm}
\begin{enumerate}
\item[{\rm (a)}] If $X \cap B^{\hspace{.2mm}e} = B^{\hspace{.2mm}e -1}$ then $\widetilde{\beta}_i(Y) =
\widetilde{\beta}_i(X)$ for all $i$. \item[{\rm (b)}] If $X \cap B^{\hspace{.2mm}e} = S^{\hspace{.2mm}e -1}$ then
\begin{enumerate}
\item[\mbox{either}] ~  $\widetilde{\beta}_i(Y) - \widetilde{\beta}_i(X) = \left\{\begin{array}{rll}
+1 & \mbox{if } & i=e \\
0 & \mbox{if } & i\neq e,
\end{array}
\right.$ \item[\mbox{or}] ~ $\widetilde{\beta}_i(Y) - \widetilde{\beta}_i(X) = \left\{\begin{array}{rll}
-1 & \mbox{if } & i=e-1 \\
0 & \mbox{if } & i\neq e-1.
\end{array}
\right.$
\end{enumerate}
\end{enumerate}
\end{lemma}

\noindent {\bf Proof.} When $e\neq 0$, this is immediate from Mayer-Vietoris theorem for reduced (simplicial)
homology. When $e=0$, the hypothesis says that $Y$ is the disjoint union of $X$ and a point, so that the result
is trivial in this case (and the first alternative holds). \hfill $\Box$

\begin{lemma} \label{L8}
Let $X$, $Y$ be simplicial complexes of dimension $d$ such that $Y$ is obtained from $X$ by a single bistellar
move $\alpha \mapsto \beta$, say of index $t$ $(0 \leq t \leq d\,)$. Then, for any set $A$, the reduced Betti
numbers of $X[A]$ and $Y[A]$ are related as follows. \vspace{-1mm}
\begin{enumerate}
\item[$(i)$] If $A \supseteq \beta$, $A \cap \alpha = \emptyset$ then \vspace{-1mm}
\begin{enumerate}
\item[\mbox{either}] ~  $\widetilde{\beta}_i(Y[A]) - \widetilde{\beta}_i(X[A]) = \left\{\begin{array}{rll}
+1 & \mbox{if } & i=t \\
0 & \mbox{if } & i\neq t,
\end{array}
\right.$ \item[\mbox{or}] ~ $\widetilde{\beta}_i(Y[A]) - \widetilde{\beta}_i(X[A]) = \left\{\begin{array}{rll}
-1 & \mbox{if } & i=t-1 \\
0 & \mbox{if } & i\neq t-1.
\end{array}
\right.$
\end{enumerate}
\vspace{-2mm} \item[$(ii)$] If $A \supseteq \alpha$, $A \cap \beta = \emptyset$ then \vspace{-2mm}
\begin{enumerate}
\item[\mbox{either}] ~  $\widetilde{\beta}_i(Y[A]) - \widetilde{\beta}_i(X[A]) = \left\{\begin{array}{rll}
-1 & \mbox{if } & i= d-t \\
0 & \mbox{if } & i\neq d-t,
\end{array}
\right.$
\item[\mbox{or}] ~ $\widetilde{\beta}_i(Y[A]) - \widetilde{\beta}_i(X[A]) = \left\{\begin{array}{rll}
+1 & \mbox{if } & i=d-t-1 \\
0 & \mbox{if } & i\neq d-t-1.
\end{array}
\right.$
\end{enumerate}
\vspace{-2mm} \item[$(iii)$] In all other cases, $\widetilde{\beta}_i(Y[A]) = \widetilde{\beta}_i(X[A])$ for all
$i$.
\end{enumerate}
\end{lemma}

\noindent {\bf Proof.} If $A \supseteq \beta$, $A \cap \alpha = \emptyset$ then we have $Y[A] = X[A] \cup
\overline{\beta}$ and $X[A] \cap \overline{\beta} = \partial \beta$. If $A \supseteq \alpha$, $A \cap \beta =
\emptyset$ then we have $X[A] = Y[A] \cup \overline{\alpha}$ and $Y[A] \cap \overline{\alpha} = \partial \alpha$.
So, the result is immediate from Lemma \ref{L7} in these cases. (Remember that $\dim(\beta) =t$ and
$\dim(\alpha)=d-t$.)

If $A \supseteq \alpha\cup\beta$ then $Y[A]$ is obtained from $X[A]$ by the bistellar move $\alpha \mapsto
\beta$. Hence $Y[A]$ and $X[A]$ are homeomorphic in this case, and the result follows.

If $A$ contains neither $\alpha$ nor $\beta$, then $Y[A] = X[A]$, and the result is trivial.

If $A \supseteq \beta$ and $\alpha_0 := A\cap \alpha$ is a proper non-empty subset of $\alpha$, then $Y[A] =
X[A]\cup B^{\hspace{.2mm}e}$ and $X[A]\cap B^{\hspace{.2mm}e} = B^{\hspace{.2mm}e -1}$, where $e=t + \#(\alpha_0)
> 0$, $B^{\hspace{.2mm}e} = \overline{\alpha_0\cup\beta}$, $B^{\hspace{.2mm}e-1} = \overline{\alpha_0} \ast
\partial \beta$. Hence the result follows from Lemma \ref{L7}.

If $A \supseteq \alpha$ and $\beta_0 := A\cap \beta$ is a proper non-empty subset of $\beta$, then $X[A] =
Y[A]\cup B^{\hspace{.2mm}e}$ and $Y[A]\cap B^{\hspace{.2mm}e} = B^{\hspace{.2mm}e-1}$, where $e=d-t + \#(\beta_0)
> 0$, $B^{\hspace{.2mm}e} = \overline{\alpha\cup\beta_0}$, $B^{\hspace{.2mm}e-1} = \overline{\beta_0} \ast
\partial \alpha$. Hence the result follows from Lemma \ref{L7} in this case also. \hfill $\Box$

\bigskip

Now, we are in a position to prove a crucial result on the sigma-vectors of $k$-stellated spheres\,:

\begin{theo}  \label{P19}
For $k\geq 1$, let $S$ be an $m$-vertex $k$-stellated triangulated sphere of dimension $d \geq 2k-1$. Then, with
respect to any field, the sigma-vector of $S$ is related to its $g$-vector by\,: \vspace{-2mm}
\begin{enumerate}
\item[{\rm (a)}] $\sigma_i =0$ \, for \, $k \leq i \leq d-k-1$, \vspace{-1mm}
\item[{\rm (b)}] ${\displaystyle \sum_{i=0}^{l}(-1)^{l-i} \sigma_i \leq \frac{m+1}{d+3}
\sum_{i=0}^{l+1}(-1)^{l+1-i} \frac{g_i}{{d+2 \choose i}}}$ \, for \, $0 \leq l \leq k-2$, \, and \vspace{-1mm}
\item[{\rm (c)}] ${\displaystyle \sum_{i=0}^{l}(-1)^{l-i} \sigma_i = \frac{m+1}{d+3} \sum_{i=0}^{l+1}(-1)^{l+1-i}
\frac{g_i}{{d +2 \choose i}}}$ \, for \, $k-1 \leq l \leq d-k-1$.
\end{enumerate}
\end{theo}

\noindent {\bf Proof.} Induction on the minimum number $\lambda(S)$ of bistellar moves (of index $<k$) needed to
obtain $S$ from $S^{\hspace{.2mm}d}_{d+2}$. If $\lambda(S) = 0$, then $S= S^{\hspace{.2mm}d}_{d+2}$. In this
case, $\sigma_i(S) = -\delta_{i0}$ for $0\leq i<d$, $g_i(S) = \delta_{i0}$ for $0\leq i\leq d+1$, and $m=d+2$.
So, the result is trivial in this case. Now, assume $\lambda(S) >0$. Then $S$ is obtained from a $k$-stellated
$d$-sphere $S^{\hspace{.2mm}\prime}$ (with $\lambda(S^{\hspace{.2mm}\prime}) = \lambda(S)-1$) by a single
bistellar move $\alpha \mapsto \beta$, say of index $t$ ($0\leq t <k$).

For $k \leq i\leq d-k-1$, we have $t < i < d-t-1$ and hence by Lemma \ref{L8}, $\widetilde{\beta}_i(S[A]) =
\widetilde{\beta}_i(S^{\hspace{.2mm}\prime}[A])$ for all subsets $A$ of $V(S)$. Taking an appropriate weighted
sum of these equalities over all sets $A$, we get $\sigma_i(S) = \sigma_i(S^{\hspace{.2mm}\prime}) = 0$ for
$k\leq i \leq d-k-1$. Here the last equality is by induction hypothesis. This proves part (a).

Fix an index $l$ such that $0\leq l\leq d-k-1$ and let $0\leq i \leq l$. Then $i \leq d-k-1 < d-t-1$, and hence
by Lemma \ref{L8}, $\widetilde{\beta}_i(S[A]) = \widetilde{\beta}_i(S^{\hspace{.2mm}\prime}[A])$ unless
$A\supseteq \beta$, $A\cap \alpha = \emptyset$. So, let ${\cal A}$ be the collection of all subsets $A$ of $V(S)$
such that $A\supseteq \beta$ and $A\cap \alpha = \emptyset$. For $A \in {\cal A}$, we have
$\widetilde{\beta}_i(S[A]) = \widetilde{\beta}_i(S^{\hspace{.2mm}\prime}[A])$ for $i\neq t, t-1$, and either
$\widetilde{\beta}_t(S[A])= \widetilde{\beta}_t(S^{\hspace{.2mm}\prime}[A]) +1$, $\widetilde{\beta}_{t-1}(S[A])=
\widetilde{\beta}_{t-1}(S^{\hspace{.2mm}\prime}[A])$, or else $\widetilde{\beta}_t(S[A])=
\widetilde{\beta}_t(S^{\hspace{.2mm}\prime}[A])$, $\widetilde{\beta}_{t-1}(S[A])= \widetilde{\beta}_{t
-1}(S^{\hspace{.2mm}\prime}[A])-1$. Let ${\cal A}^{\hspace{.2mm}+}$ be the set of all $A\in {\cal A}$ for which
the first alternative holds. Then we get\,:
\begin{eqnarray} \label{eq10}
\sum_{i=0}^l (-1)^{l-i}\left[\,\widetilde{\beta}_i(S[A]) - \widetilde{\beta}_i(S^{\hspace{.2mm}\prime}[A])\right]
= \left\{\begin{array}{cl}
0 & \mbox{if } ~ A \not\in {\cal A} ~ \mbox{or if } ~ l < t-1 \\
0 & \mbox{if } ~ A \in {\cal A}^{\hspace{.2mm}+} \, \mbox{ and  } ~ l = t-1 \\ (-1)^{l-t} & \mbox{otherwise.}
\end{array}
\right.
\end{eqnarray}

First consider the case $l\leq t-1$  (which can occur only for $l < k-1$ and $t>0$). Then (\ref{eq10}) implies
$$
\sum_{i=0}^l (-1)^{l-i} \widetilde{\beta}_i(S[A]) \leq \sum_{i=0}^l (-1)^{l-i}
\widetilde{\beta}_i(S^{\hspace{.2mm}\prime}[A]) ~~ \mbox{for all} ~~ A \subseteq V(S).
$$
Taking the appropriate weighted sum of these inequalities over all $A$, we get
\begin{eqnarray*}
\sum_{i=0}^l (-1)^{l-i}\sigma_i(S) ~ \leq ~ \sum_{i=0}^l (-1)^{l-i} \sigma_i(S^{\hspace{.2mm}\prime}) &\leq
&\frac{m+1}{d+3} \sum_{i=0}^{l+1}(-1)^{l+1-i}\,
\frac{g_i(S^{\hspace{.2mm}\prime})}{{d+2 \choose i}} \\
&= &\frac{m+1}{d+3} \sum_{i=0}^{l+1}(-1)^{l+1-i}\, \frac{g_i(S)}{{d+2 \choose i}},
\end{eqnarray*}
where the second inequality is by induction hypothesis and the final equality holds since by Lemma \ref{L3}, we
have $g_i(S) = g_i(S^{\hspace{.2mm}\prime})$ for $i\leq l+1 \leq t$. This completes the induction step in this
case.

Next consider the case $t \leq l \leq d-k-1$. In this case, (\ref{eq10}) says\,:
$$
\sum_{i=0}^l (-1)^{l-i} \left[\,\widetilde{\beta}_i(S[A]) -
\widetilde{\beta}_i(S^{\hspace{.2mm}\prime}[A])\right] = \left\{\begin{array}{cl}
(-1)^{l-t} & \mbox{if } ~ A \in {\cal A} \\
0 & \mbox{otherwise.}
\end{array}
\right.
$$
Notice that there are exactly ${m-d-2 \choose j-t-1}$ $j$-sets in ${\cal A}$. Therefore, adding these equations
over all $j$-subsets of $V := V(S)$ we get (for $0\leq j \leq m$)\,:
\begin{eqnarray} \label{eq11}
\sum_{i=0}^l (-1)^{l-i} \sum_{A\in{V \choose j}} \left[\,\widetilde{\beta}_i(S[A]) -
\widetilde{\beta}_i(S^{\hspace{.2mm}\prime}[A])\right] = (-1)^{l-t}{m-d-2 \choose j-t-1}.
\end{eqnarray}
First suppose $t>0$, so that $V(S^{\hspace{.2mm}\prime}) = V(S)$ has size $m$. Dividing Equation (\ref{eq11}) by
${m\choose j}$ and adding over all $j$, we get\,:
$$
\sum_{i=0}^l (-1)^{l-i} \left[\sigma_i(S) - \sigma_i(S^{\hspace{.2mm}\prime})\right] =
(-1)^{l-t}\sum_{j=0}^m\frac{{m-d-2 \choose j-t-1}}{{m\choose j}} = (-1)^{l-t}\sum_{i=0}^{m-d-2} \frac{{m-d-2
\choose i}}{{m \choose t+1+i}}.
$$
Now, the binomial coefficients satisfy the following well known identity. For any three non-negative integers $p,
q, r$, we have
\begin{eqnarray} \label{eq12}
\sum_{i=0}^{p} \frac{{p \choose i}}{{p+q+r \choose r+i}} = \frac{p+q+r+1}{q+r+1}\times \frac{1}{{q+r \choose r}}.
\end{eqnarray}
Substituting $m-d-2$, $d+1-t$ and $t+1$ for $p$, $q$, $r$ respectively, we get
$$
\sum_{i=0}^l (-1)^{l-i} \left[\sigma_i(S) - \sigma_i(S^{\hspace{.2mm}\prime})\right] = (-1)^{l-t}\,
\frac{m+1}{d+3}\frac{1}{{d+2 \choose t+1}}.
$$
Therefore, induction hypothesis gives the following inequality for $0<t\leq l\leq d-k-1$ (with equality for
$k-1\leq l\leq d-k-1$)\,:
\begin{eqnarray*}
\sum_{i=0}^l (-1)^{l-i}\sigma_i(S) &\leq & (-1)^{l-t} \, \frac{m+1}{d+3}\frac{1}{{d+2 \choose t+1}} +
\frac{m+1}{d+3} \sum_{i=0}^{l+1} (-1)^{l+1-i} \,
\frac{g_i(S^{\hspace{.2mm}\prime})}{{d+2 \choose i}} \\
&= &\frac{m+1}{d+3} \sum_{i=0}^{l+1}(-1)^{l+1-i}\, \frac{g_i(S^{\hspace{.2mm}\prime}) + \delta_{i, t+1}}{{d+2
\choose
i}}  \\
&= & \frac{m+1}{d+3} \sum_{i=0}^{l+1}(-1)^{l+1-i} \, \frac{g_i(S)}{{d+2 \choose i}},
\end{eqnarray*}
where the last equality holds since by Lemma \ref{L3}, we have $g_i(S) = g_i(S^{\hspace{.2mm}\prime}) +
\delta_{i, t+1}$ as $d\geq 2k-1$ (so that $d\neq 2t$) and $i-1 \leq l \leq d-k < d-t$. This completes the
induction step in the second case.

Finally, consider the case $0=t \leq l\leq d-k-1$. In this case, $\beta$ is a vertex of $S$ not in
$S^{\hspace{.2mm}\prime}$. Let $V^{\hspace{.2mm}\prime} = V\setminus \{\beta\}$ be the vertex set of
$S^{\hspace{.2mm}\prime}$. (Thus, $S^{\hspace{.2mm}\prime}$ has $m-1$ vertices in this case.) Dividing Equation
(\ref{eq11}) (with $t=0$) by ${m \choose j}$ and adding over all $j$ ($0\leq j\leq m$) we get (in view of
(\ref{eq12}))

\begin{eqnarray*}
\sum_{i=0}^l (-1)^{l-i}\sigma_i(S) = (-1)^{l} \, \frac{m+ 1}{(d+3)(d+2)} + \sum_{i=0}^l (-1)^{l-i} \sum_{j=0}^m
\frac{1}{{m \choose j}} \sum_{A\in{V \choose j}} \widetilde{\beta}_i(S^{\hspace{.2mm}\prime}[A]).
\end{eqnarray*}
Now, for each $i\leq l$,
\begin{eqnarray*}
\sum_{j=0}^m \frac{1}{{m \choose j}} \sum_{A\in{V \choose j}} \widetilde{\beta}_i(S^{\hspace{.2mm}\prime}[A]) &=
& \sum_{j=0}^{m-1} \frac{1}{{m \choose j}} \sum_{A\in{V \choose j}}
\widetilde{\beta}_i(S^{\hspace{.2mm}\prime}[A]) \, = \, \sum_{j=0}^{m-1} \frac{1}{{m \choose j}} \sum_{B\in
{V^{\hspace{.1mm}\prime} \choose j} \sqcup {V^{\hspace{.1mm}\prime} \choose j-1}}
\hspace{-2mm}\widetilde{\beta}_i(S^{\hspace{.2mm}\prime}[B]) \\
&= & \sum_{j=0}^{m-1} \left[\frac{1}{{m \choose j}} + \frac{1}{{m \choose j+1}}\right] \sum_{B\in
{V^{\hspace{.1mm}\prime} \choose j}} \hspace{-2mm} \widetilde{\beta}_i(S^{\hspace{.2mm}\prime}[B]) \, = \,
\frac{m+1}{m} \, \sigma_i(S^{\hspace{.2mm}\prime}).
\end{eqnarray*}
Since $i\leq l\leq d-k-1 <d$, $\widetilde{\beta}_i(S^{\hspace{.2mm} \prime}) = 0$. This justifies the first and
third equalities above. The second equality holds since, by definition, $S^{\hspace{.2mm} \prime}[A] =
S^{\hspace{.2mm}\prime}[B]$, where $B = A \cap V^{\hspace{.1mm} \prime}$, and the map $A\mapsto A\cap
V^{\hspace{.1mm} \prime}$ is a bijection between ${V \choose j}$ and ${\,V^{\hspace{.1mm}\prime} \choose j}
\sqcup {V^{\hspace{.1mm}\prime} \choose j-1\,}$. The last equality is because of the trivial identity
$\frac{1}{{m \choose j}} + \frac{1}{{m \choose j+1}} = \frac{m+1}{m} \frac{1}{{m -1 \choose j}}$. So, we have
(when $t=0$)\,:
$$
\sum_{i=0}^l (-1)^{l-i}\sigma_i(S) = (-1)^{l} \, \frac{m+ 1}{(d+3)(d+2)} + \frac{m+1}{m}\sum_{i=0}^l (-1)^{l-i}
\sigma_i(S^{\hspace{.2mm}\prime}), ~ 0\leq l\leq d-k-1.
$$
Now, since $S^{\hspace{.2mm}\prime}$ has $m-1$ vertices, induction hypothesis gives
$$
\sum_{i=0}^l (-1)^{l-i} \sigma_i(S^{\hspace{.2mm}\prime}) \leq \frac{m}{d+3}\sum_{i=0}^{l+1} (-1)^{l+1-i} \,
\frac{g_i(S^{\hspace{.2mm}\prime})}{{d+2 \choose i}}\,,
$$
with equality for $k-1 \leq l \leq d-k-1$. Therefore, we get
\begin{eqnarray*}
\sum_{i=0}^l (-1)^{l-i}\sigma_i(S) & \leq & (-1)^{l} \, \frac{m+ 1}{(d+3)(d+2)} + \frac{m+1}{d+3}\sum_{i=0}^{l+1}
(-1)^{l+1-i} \,
\frac{g_i(S^{\hspace{.2mm}\prime})}{{d+2 \choose i}} \\
&=& \frac{m+1}{d+3}\sum_{i=0}^{l+1} (-1)^{l+1-i} \,
\frac{g_i(S^{\hspace{.2mm}\prime})+ \delta_{i1}}{{d+2 \choose i}}  \\
&=& \frac{m+1}{d+3}\sum_{i=0}^{l+1} (-1)^{l+1-i} \, \frac{g_i(S)}{{d+2 \choose i}}\,,
\end{eqnarray*}
with equality for $k-1 \leq l\leq d-k-1$. This completes the induction in the last case, thus proving (b) and
(c). \hfill $\Box$

\bigskip

Now, the following key result on the mu-vectors of 2-neighbourly members of ${\cal W}_k(d)$ is more or less
immediate.

\begin{theo}  \label{P20}
Let $M \in {\cal W}_k(d)$ be $2$-neighbourly with $d\geq 2k\geq 2$. Then the mu-vector of $M$ is related to its
$g$-vector as follows\,: \vspace{-1mm}
\begin{enumerate}
\item[{\rm (a)}] $\mu_i =0$ \, for \, $k+1 \leq i \leq d-k-1$, \vspace{-1mm}
\item[{\rm (b)}]  ${\displaystyle \sum_{i=1}^{l}(-1)^{l-i}\mu_i \leq \frac{g_{l+1}}{{d+2 \choose l+1}}}$ \, for
\, $1 \leq l \leq k-1$, \, and  \vspace{-2mm}
\item[{\rm (c)}] ${\displaystyle \sum_{i=1}^{l}(-1)^{l-i} \mu_i = \frac{g_{l+1}}{{d+2 \choose l+1}}}$, \,  for \,
$k \leq l \leq d-k-1$.
\end{enumerate}
\end{theo}

\noindent {\bf Proof.} Let $m$ be the number of vertices of $M$ and let $L_x$ be the link of $x$ in $M$ for $x\in
V(M)$. Then each $L_x$ is a $k$-stellated sphere, of dimension $d-1\geq 2k-1$, on exactly $m-1$ vertices.
Therefore, for $k+1 \leq i \leq d-k-1$, $\sigma_{i-1}(L_x) = 0$ by Theorem \ref{P19}\,(a). Taking the sum of
these equations over all $x\in V(M)$, we get $\mu_i(M) = 0$ for $k+1 \leq i \leq d-k-1$. This proves part (a).

Also, by Theorem \ref{P19}\,(b) and (c), we have, for $1\leq l\leq d-k-1$,
$$
\sum_{i=1}^l (-1)^{l-i} \, \sigma_{i-1}(L_x) \leq \frac{m}{d+2} \sum_{i=0}^l (-1)^{l-i} \, \frac{g_i(L_x)}{{d+1
\choose i}},
$$
with equality for $k\leq l \leq d-k-1$. Adding these over all $x\in V(M)$, and dividing the result by $m$, we get
(in view of Lemma \ref{L4})\,:
$$
\sum_{i=1}^l (-1)^{l-i} \, (\mu_{i} - \delta_{i1}) \leq \frac{1}{d+2} \sum_{i=0}^l \frac{(-1)^{l-i}}{{d+1 \choose
i}}((d+2-i)g_i + (i+1)g_{i+1}),
$$
with equality for $k\leq l \leq d-k-1$. That is,
$$
\sum_{i=1}^l (-1)^{l-i} \, \mu_{i} \leq (-1)^{l-1} + \sum_{i=0}^l (-1)^{l-i} \left(\frac{g_i}{{d+2 \choose i}} +
\frac{g_{i+1}}{{d+2 \choose i+1}}\right) = \frac{g_{l+1}}{{d+2 \choose l+1}} ~ (\mbox{since } ~ g_0 =1),
$$
with equality for $k\leq l \leq d-k-1$. This proves (b) and (c). \hfill $\Box$

\bigskip

Now, we can prove one of the main results of this paper.

\begin{theo}[A lower bound theorem for ${\cal W}_k(d)$]  \label{P21}
Let $M \in {\cal W}_k(d)$ be  $2$-\hspace{-1.5mm} neighbourly. Then the $g$-vector of $M$ is related to its Betti
numbers  as follows\,: \vspace{-3mm}
\begin{enumerate}
\item[{\rm (a)}] if \, $d = 2k$,  \, then \, $g_{l+1} \geq {d+2 \choose l+1}  {\displaystyle
\sum_{i=1}^{l}(-1)^{l-i}\beta_i }$ \, for \, $1 \leq l \leq k-1$,
 \vspace{-3mm}
\item[{\rm (b)}] if \, $d \geq 2k+1$, \, then \,  $g_{l+1} \geq {d+2 \choose l+1}  {\displaystyle
\sum_{i=1}^{l}(-1)^{l-i}\beta_i }$ \, for \, $1 \leq l \leq k$,
 \vspace{-3mm}
\item[{\rm (c)}] if \, $d \geq 2k+2$, \, then \, $g_{l+1} = {d+2 \choose l+1}  {\displaystyle
\sum_{i=1}^{l}(-1)^{l-i}\beta_i }$ \, for \, $k \leq l \leq d-k-1$, \, and
\item[{\rm (d)}] if \, $d\geq 2k+2$, \, then \, $\beta_{i} = 0$ \, for \, $k+1 \leq i \leq d-k-1$.
\end{enumerate}
\end{theo}

\noindent {\bf Proof.} By Theorems \ref{P16}\,(b) and \ref{P20}\,(a), we have $0\leq \beta_i \leq \mu_i=0$ for
$k+1 \leq i \leq d-k-1$. This proves part (d).

Since $M$ is 2-neighbourly, it is connected. Therefore, $\beta_0 = 1 = \mu_0$. Hence, Theorem \ref{P16}\,(a)
yields $\sum_{i= 1}^{l}(- 1)^{l-i} \beta_i \leq  { \sum_{i=1}^{l}(-1)^{l-i} \mu_i }$ for $l \geq 1$. Therefore,
parts (a) and (b) are immediate from parts (b) and (c) of Theorem \ref{P20}.

If $k < l \leq d-k-1$, then $\beta_l = \mu_l$ by part (d), and hence by Theorem \ref{P16}\,(d), $H_l(Y) \to
H_l(M)$ is injective for any induced subcomplex $Y$ of $M$. Hence, by Theorem \ref{P16}\,(c), ${ \sum_{i=
1}^{l}(- 1)^{l-i} \beta_i } = {\sum_{i= 1}^{l}(-1)^{l-i}\mu_i }$. But, ${\sum_{i= 1}^{l}(- 1)^{l-i}\mu_i } =
g_{l+1}/{d+2 \choose l+1}$ by Theorem \ref{P20}\,(c). This proves (c) in case $k<l\leq d-k-1$. Finally, when $k+1
\leq d-k-1$, we have $\beta_{k+1} =0= \mu_{k+1}$ by part (d) and Theorem \ref{P20}\,(a), hence $H_k(Y) \to
H_k(M)$ is injective for any induced subcomplex $Y$ of $M$ (by Theorem \ref{P16}\,(d)). Therefore, by Theorems
\ref{P16}\,(c) and \ref{P20}\,(c), ${ \sum_{i=1}^{k}(- 1)^{k-i} \beta_i } = {\sum_{i= 1}^{k}(-1)^{k -i}\mu_i } =
g_{k+1}/{d+2 \choose k+1}$, completing the proof of part (c).  \hfill $\Box$

\bigskip

We also need the following elementary result.

\begin{lemma} \label{L9}
Let $X$ be an $(l+1)$-neighbourly simplicial complex. Then the beta- and mu-vectors of $X$ satisfy $\beta_i = 0 =
\mu_i$ for $1\leq i \leq l-1$.
\end{lemma}

\noindent {\bf Proof.} The $l$-skeleton of $X$ agrees with that of the standard ball of dimension $f_0(X)-1$.
Since the ball is homologically trivial and the $i^{\rm th}$ homology of a simplicial complex is the same as the
$i^{\rm th}$ homology of its $(i+1)$-skeleton, it follows that $\beta_i(X) = 0$ for $1\leq i\leq l-1$.

Also, for any vertex $x$ of $X$, the link $L_x$ of $x$ in $X$ is $l$-neighbourly. Therefore, by the same
argument, we have, for $1\leq i \leq l-1$, $\widetilde{\beta}_{i-1}$ of any induced subcomplex of $L_x$ is $=0$,
except that $\widetilde{\beta}_{0} = -1$ for the empty subcomplex. Therefore, taking an appropriate weighted sum,
we get $\sigma_{i-1}(L_x) = - \delta_{i1}$ for $1\leq i \leq l-1$ and $x\in V(X)$. Adding over all $x\in V(X)$,
we get $\mu_i(X) =0$ for $1\leq i \leq l-1$. \hfill $\Box$

\bigskip

The following result (which is the first known combinatorial criterion for tightness) is due to K\"{u}hnel
\cite{k95}.

\begin{lemma}[K\"{u}hnel] \label{L10}
For $k\geq 1$, let $M$ be a $(k+1)$-neighbourly triangulation of an $\mathbb{F}$-orientable closed manifold of
dimension $2k$. Then $M$ is $\mathbb{F}$-tight.
\end{lemma}

\noindent {\bf Proof.} Since $M$ is at least 2-neighbourly, it is connected. Therefore, $\mu_0 = 1 = \beta_0$. By
Lemma \ref{L9}, $\mu_i = 0=\beta_i$ for $1\leq i\leq k-1$. So, by duality (Theorem \ref{P16}\,(e)), $\mu_i = 0 =
\beta_i$ for $k+1 \leq i \leq 2k-1$ and $\mu_{2k} = 1 = \beta_{2k}$. Thus $\mu_i = \beta_i$ for all $i$, except
possibly for $i=k$. But then, the equality $\sum_{i= 0}^{2k} (-1)^i \mu_i = \sum_{i=0}^{2k} (-1)^i \beta_i$ from
Theorem \ref{P16}\,(a) implies $\mu_i = \beta_i$ for $i=k$ as well. Therefore, by Theorem \ref{P18}, $M$ is
$\mathbb{F}$-tight. \hfill $\Box$

\bigskip

The $k=1$ case of Theorem \ref{P25}\,(a) below is essentially due to Effenberger \cite{ef}. This paper was
largely motivated by a desire to understand and generalise Effenberger's result.

\begin{theo}[A combinatorial criterion for tightness] \label{P25}
Let $M \in {\cal W}_k^{\ast}(d)$ be $\mathbb{F}$-orientable. Then we have\,: \vspace{-2mm}
\begin{enumerate}
\item[{\rm (a)}]  if $d \neq 2k+1$ then $M$ is $\mathbb{F}$-tight, and \vspace{-2mm}
\item[{\rm (b)}] if $d = 2k+1$, then $M$ is $\mathbb{F}$-tight if and only if \, $\beta_k(M; \mathbb{F}) =
\frac{{n-k-3 \choose k+1}}{{2k+3 \choose k+1}}$, where $n = f_0(M)$.
\end{enumerate}
\end{theo}
(Note that, when $k\geq 2$, $M$ is by assumption at least 3-neighbourly, and hence simply connected. Thus, the
hypothesis of orientability is automatic for $k\geq 2$.)

\bigskip

\noindent {\bf Proof.} This is trivial if $M = S^{\hspace{.2mm}d}_{d +2}$. So, assume $M \neq
S^{\hspace{.2mm}d}_{d+2}$. Then $d\geq 2k$. If $d=2k$, then the result follows from Lemma \ref{L10}. So, assume
that $d\geq 2k+1$.

Since $M$ is $\mathbb{F}$-orientable, the duality result of Theorem \ref{P16}\,(e) applies. Since $M$ is
connected, we have $\beta_0 = 1 = \mu_0$ and hence $\beta_d = 1 = \mu_d$. By Lemma \ref{L9}, $\beta_i = 0 =
\mu_i$ for $1\leq i \leq k-1$, hence also for $d-k+1\leq i \leq d-1$. Since $M$ is $(k+1)$-neighbourly, we have
$g_{k+1}(X) = {n+k-d-2 \choose k+1}$, where $n=f_0(M)$. We also have $\beta_k = g_{k+1}/{d+2 \choose k+1} =
\mu_k$ (by hypothesis and Theorem \ref{P20}\,(c) when $d = 2k+1$; by Theorems \ref{P20}\,(c) and \ref{P21}\,(c)
when $d\geq 2k+2$). Hence, $\beta_{d-k} = g_{k+1}/{d+2 \choose k+1} = \mu_{d-k}$. By Theorems \ref{P20}\,(a) and
\ref{P21}\,(d), we also have $\beta_i = 0 = \mu_i$ for $k+1 \leq i \leq d-k-1$. Thus, $\beta_i = \mu_i$ for all
$i$. Hence $M$ is $\mathbb{F}$-tight by Theorem \ref{P18} (and, when $d\neq 2k+1$, this argument applies to all
fields $\mathbb{F}$).

For the converse statement in part (b), note that - more generally - for any $M \in {\cal W}_k^{\ast}(d)$ with $d
\geq 2k+1$, Lemma \ref{L9} and Theorem \ref{P20}\,(c) imply that $\mu_k = g_{k+ 1}/{d+2 \choose k+1}$. Therefore,
for $M$ to be $\mathbb{F}$-tight, we must have (by Theorem \ref{P18}) $\beta_k = g_{k+1}/{d+2 \choose k+1}$ as
well. Thus, for $d=2k+1$, $\beta_k = {n+k-d-2 \choose k+1}/{d+2 \choose k+1} = {n-k-3 \choose k+1}/{2k+3 \choose
k+1}$. \hfill $\Box$


\bigskip

The `$d=4$' case of the following theorem is due to Walkup \cite{wa} and K\"{u}hnel \cite{k95} (cf.
\cite[Proposition 2]{bd10}). Part (b) of this theorem is due to Lutz, Sulanke and Swartz \cite{lss}.

\begin{theo}[A lower bound theorem for triangulated manifolds]  \label{P23}
Let $M$ be a connected closed triangulated manifold of dimension $d\geq 3$. Let $\beta_1 = \beta_1(M;
\mathbb{Z}_2)$. Then the face vector of $M$ satisfies\,:
\begin{enumerate}
\item[{\rm (a)}] $f_j \geq \left\{\begin{array}{ll} {\,d+1\, \choose j}f_0 +j{\,d+2\, \choose j+1}(\beta_1 -1),
& \mbox{if } ~ 1 \leq j < d,  \\[2mm]
df_0 + (d-1)(d+2)(\beta_1-1), & \mbox{if } ~ j = d.
\end{array}
\right. $ \item[{\rm (b)}] ${\,f_0-d-1\, \choose 2} \geq {\,d+2\, \choose 2}\beta_1$.
\end{enumerate}
When $d\geq 4$, equality holds in {\rm (a)} $($for some $j \geq 1$, equivalently, for all $j\,)$ if and only if
$M\in {\cal W}_1(d)$, and equality holds in {\rm (b)} if and only if $M$ is a $2$-neighbourly member of ${\cal
W}_1(d)$.
\end{theo}

\noindent {\bf Proof.} The face vector of any 1-stellated $(d-1)$-sphere $S$ satisfies (cf. \cite{bd13v2}, for
instance)
$$
f_{j-1}(S) = \left\{\begin{array}{ll} {d \choose j-1\,}f_0(S) -(j-1){\,d+1\, \choose j}, & \mbox{if }
~ 1 \leq j < d,  \\[2mm]
(d-1)f_0(S) - (d-2)(d+1), & \mbox{if } ~ j = d.
\end{array}
\right.
$$
The lower bound theorem of Kalai \cite{ka} says that if $L$ is any triangulated closed manifold of dimension
$d-1$ with $f_0(L) = f_0(S)$, then $f_{j-1}(L) \geq f_{j-1}(S)$ for $1\leq j\leq d$. Also, when $d\geq 4$,
equality holds here (for some $j>1$, equivalently, for all $j$) if and only if $L$ is an 1-stellated sphere.
Applying this result to the vertex links $L_x$, $x \in V(M)$, of $M$, we get
$$
f_{j-1}(L_x) \geq \left\{\begin{array}{ll} {d \choose \,j-1\,}f_0(L_x) -(j-1){\,d+1\, \choose j},
& \mbox{if } ~ 1 \leq j < d,  \\[2mm]
(d-1)f_0(L_x) - (d-2)(d+1), & \mbox{if } ~ j = d.
\end{array}
\right.
$$
But $\sum_{x\in V(M)}f_{j-1}(L_x) = (j+1)f_j$, $1\leq j \leq d$. Therefore, adding the above inequalities over
all vertices $x$ of $M$, we get that the face vector of $M$ satisfies\,:

\begin{eqnarray} \label{eq13}
f_{j} \geq \left\{\begin{array}{ll} \frac{2}{j+1}{d \choose \,j-1\,}f_1 -\frac{j-1}{j+1}{\,d+1\,
\choose j}f_0, & \mbox{if } ~ 1 \leq j < d,  \\[2mm]
\frac{2d-2}{d+1}f_1 - (d-2)f_0, & \mbox{if } ~ j = d.
\end{array}
\right.
\end{eqnarray}
Also, when $d\geq 4$, equality holds in (\ref{eq13}) for some $j> 1$ (equivalently, for all $j$) if and only if
all the vertex links of $M$ are 1-stellated, i.e., if and only if $M\in {\cal W}_1(d)$.

Now, Theorem 5.2 of Novik and Swartz \cite{ns} says that $g_2(M) \geq {d+2 \choose 2}\beta_1$, i.e.,
\begin{eqnarray} \label{eq14}
f_1 \geq (d+1)f_0 + {d+2 \choose 2}(\beta_1 -1).
\end{eqnarray}
Also, when $d\geq 4$, equality holds in (\ref{eq14}) if and only if $M\in {\cal W}_1(d)$. Notice that
(\ref{eq14}) is just the case $j=1$ of part (a). Now, combining (\ref{eq13}) and (\ref{eq14}), we get all cases
of part (a), after a little simplification.

We also have $ (d+1)f_0 + {d+2\, \choose 2}(\beta_1-1) \leq f_1 \leq \frac{f_0(f_0-1)}{2}\,, $ where the first
inequality is from (\ref{eq14}) (with equality for $d\geq 4$ if and only if $M \in {\cal W}_1(d)$) and the second
inequality is trivial (with equality if and only if $M$ is 2-neighbourly). Hence we have $(d+1)f_0 + {d+2\,
\choose 2}(\beta_1-1) \leq \frac{f_0(f_0 -1)}{2}$, which simplifies to the inequality in part (b). \hfill $\Box$

\bigskip

In \cite{kl}, K\"{u}hnel and Lutz conjectured that all $\mathbb{F}$-tight triangulated  manifolds are strongly
minimal. Our next result is a powerful evidence in favour of this conjecture. It also generalises a result of
Swartz \cite[Theorem 4.7]{s09}, who proved that the tight triangulations $K^{\,d}_{2d+3}\in {\cal W}_1(d)$ (cf.
Example \ref{E3}\,(b) below) are strongly minimal for all $d\geq 2$.

\begin{cor} \label{C11}
Every $\mathbb{F}$-tight member of\,  ${\cal W}_1(d)$ is strongly minimal.
\end{cor}

\noindent {\bf Proof.} Let $M_0 \in {\cal W}_1(d)$ be $\mathbb{F}$-tight. Theorem \ref{P25} implies that $M_0$ is
$\mathbb{Z}_2$-tight. By the same theorem, $M_0$ is 2-neighbourly. Now, for $d\leq 2$, any 2-neighbourly closed
$d$-manifold is strongly minimal. (This is entirely trivial for $d=1$ since $S^1_3$ is the only 2-neighbourly
closed manifold in this case. It is only a little less trivial for $d=2$\,: the face vector is determined by
$f_0$ and $\beta_1$ in this case.) So, assume $d\geq 3$.

We claim that $M_0$ attains all the bounds in Theorem \ref{P23}. This is immediate from the theorem itself when
$d\geq 4$. If $d=3$, then - as $M_0\in {\cal W}_1(3)$ is $\mathbb{Z}_2$-tight, we have $g_2(M_0) = 10\beta_1$ by
Theorem \ref{P25}. Hence, following the proof of Theorem \ref{P23}, one sees that $M_0$ attains the bounds in
Theorem \ref{P23} in this case also.

Now, let $M$ be any triangulation of $|M_0|$. Since the Betti numbers are topological invariants, we have
$\beta_1(M; \mathbb{Z}_2) = \beta_1(M_0; \mathbb{Z}_2) = \beta_1$ (say). By Theorem \ref{P23} and the above
claim,
$$
{f_0(M) -d-1 \choose 2} \geq {d+2 \choose 2}\beta_1 = {f_0(M_0) -d-1 \choose 2}.
$$
Since, trivially, $f_0(M)$, $f_0(M_0) \geq d+2$, this implies $f_0(M) \geq f_0(M_0)$. Therefore, we get\,: $
f_j(M) \geq a_jf_0(M) + b_j \geq a_jf_0(M_0) + b_j = f_j(M_0) ~ \mbox{ for } ~ 0\leq j\leq d, $ where $a_j >0$
and $b_j$ are constants (depending only on $d$, $j$ and $\beta_1$) given by Theorem \ref{P23}. \hfill $\Box$

\bigskip

The members of ${\cal W}^{\ast}_1(d)$ were called ``{\em tight neighbourly}\," by Lutz, Sulanke and Swartz
\cite{lss}. By Theorem \ref{P25}, tight neighbourly manifolds of dimension $\neq 3$ are $\mathbb{Z}_2$-tight if
non-orientable, and tight if orientable. Part (b) of Theorem \ref{P25} gives a criterion for the tightness of a
tight neighbourly 3-manifold in terms of its first Betti number. Corollary \ref{C11} shows that tight neighbourly
triangulations of dimension $\geq 4$ are strongly minimal.

By a result of Kalai (cf. \cite[Corollary 8.4]{ka} or \cite[Proposition 3]{bd10}), for $d\geq 4$, any member of
${\cal W}_1(d)$ triangulates the connected sum of finitely many copies of $S^{1}\times S^{\hspace{.2mm}d-1}$ or
of $\TPSSD$. The following result may be compared with Kalai's theorem.

\begin{theo}  \label{P24}
Let $M \in {\cal W}_k^{\ast}(d)$, where $k\geq 2$ and $d\geq 2k+2$. Suppose $M$ is not the standard $d$-sphere.
Then $M$ has the same integral homology group as the connected sum of $\beta$ copies of $S^{\hspace{.2mm}k}\times
S^{\hspace{.2mm}d-k}$, where the positive integer $\beta$ is given in terms of the number $m$ of vertices of $M$
by the formula $\beta = {m+k-d-2 \choose k+1}/{d+2 \choose k+1}$. In consequence, we must have $m \geq 2d+4-k$
and ${d+2 \choose k+1}$ divides ${m+k-d-2 \choose k+1}$.
\end{theo}
(As to the inequality $m\geq 2d+4-k$, note that by a result of Brehm and K\"{u}hnel \cite{bk}, this lower bound
on the number of vertices holds, more generally, for any triangulation of a closed $d$-manifold which is not
$k$-connected.)

\bigskip

\noindent {\bf Proof.} Since $M$ is at least 3-neighbourly, it is simply connected and hence orientable.
Therefore, by Poincar\'e duality, the Betti numbers of $M$ with respect to any field $\mathbb{F}$ satisfy
$\beta_{d-i} = \beta_i$, $0\leq i \leq d$. Since $M$ is connected, we have $\beta_0 =1$ and hence $\beta_d =1$.
By Lemma \ref{L9}, $\beta_i =0$ for $1\leq i \leq k-1$. Hence, by duality, $\beta_i=0$ for $d-k+1 \leq i \leq
d-1$. By Theorem \ref{P21} and duality, $\beta_{d-k} = \beta_k = g_{k+1}/{d+2 \choose k+1}$ and $\beta_i =0$ for
$k+1 \leq i \leq d-k-1$. Thus, all the Betti numbers of $M$ are independent of the choice of the field
$\mathbb{F}$. Therefore, by the universal coefficients theorem, the $\mathbb{Z}$-homologies of $M$ are torsion
free, and the $\mathbb{Z}$-Betti numbers are given by the same formulae as above. Since $M$ is not the standard
$d$-sphere, and $M$ is $(k+1)$-neighbourly, we have $m\geq d+3$ and hence $g_{k+1} = {m+k-d-2 \choose k+1} > 0$.
Therefore, the value $\beta = g_{k+1}/{d+2 \choose k+1} = {m+k- d-2 \choose k+1}/{d+2 \choose k+1}$ of $\beta_k$
must be a strictly positive integer. Hence ${m+k-d-2 \choose k+1} \geq {d+2 \choose k+1}$, so that $m \geq
2d+4-k$. \hfill $\Box$

\begin{eg}[Tight triangulations of closed manifolds] \label{E3}
{\rm We have noted that $S^{\hspace{.2mm}d}_{d +2}$ is the only tight triangulation of $S^{\hspace{.2mm}d}$. This
trivial series apart, we know the following examples of tight triangulations (compare \cite{kl}). (Recall that we
write `tight' for $\mathbb{Q}$-tight ($\equiv$ $\mathbb{F}$-tight for all fields $\mathbb{F}$).)
\begin{enumerate}
\item[{\bf (a)}] By Lemma \ref{L10}, all 2-neighbourly triangulated closed 2-manifolds are tight when orientable
and $\mathbb{Z}_2$-tight when non-orientable. For $n \geq 4$, there exist $n$-vertex 2-neighbourly orientable
(respectively, non-orientable) triangulated closed 2-manifolds if and only if $n\equiv 0, 3, 4$ or 7 (mod 12)
(respectively, $n \equiv 0$ or 1 (mod 3), except for $n = 4, 7$) (cf. \cite{r}).

\item[{\bf (b)}] For each $d\geq 2$, there is a $(2d+3)$-vertex member $K^{d}_{2d+3}$ of ${\cal W}_1^{\ast}(d)$
found by K\"{u}hnel \cite{k86}. For $d\geq 3$, it is the unique non-simply connected $d$-manifold on $2d+3$
vertices (cf. \cite{bd8, css}). It is orientable (triangulates $S^{\,d-1}\times S^{1}$) for $d$ even and
non-orientable (triangulates $\TPSSDS$) for $d$ odd. By Theorem \ref{P25}, $K^{d}_{2d+3}$ is tight for $d$ even,
and $\mathbb{Z}_2$-tight for $d$ odd.

\item[{\bf (c)}] For each $d\geq 2$, there are two $(d^{\hspace{.1mm}2}+5d+5)$-vertex members
$M^{d}_{d^{\hspace{.1mm}2}+5d+5}$ and $N^{d}_{d^{\hspace{.1mm}2}+5d+5}$ of ${\cal W}_1^{\ast}(d)$ found by Datta
and Singh \cite{ds2}. They are orientable (triangulate $(S^{d-1}\times S^1)^{\#d^2+5d+6}$) for $d$ even and
non-orientable (triangulate $(\TPSSDS)^{\#d^2+5d+6}$) for $d$ odd. By Theorem \ref{P25}, they are tight for $d$
even, and $\mathbb{Z}_2$-tight for $d$ odd.

\item[{\bf (d)}] \begin{enumerate} \item[(i)] The 15-vertex triangulation of $(\TPSSF)^{\#3}$ obtained in
\cite{bd10} is in ${\cal W}_1^{\ast}(4)$. In \cite{si}, Nitin Singh, a student of the second author, modified
this construction to obtain another 15-vertex triangulation of $(\TPSSF)^{\#3}$ in ${\cal W}_1^{\ast}(4)$. Both
are $\mathbb{Z}_2$-tight by Theorem \ref{P25}.
\item[(ii)] Also, Singh found \cite{si} ten 15-vertex triangulations of $(S^{\,3} \times S^{1})^{\#3}$ in ${\cal
W}_1^{\ast}(4)$. All are tight by Theorem \ref{P25}.

\item[(iii)] Recently, Datta and Singh found \cite{ds} a 21-vertex triangulation of $(S^{\,3} \times
S^{1})^{\#8}$ in ${\cal W}_1^{\ast}(4)$. It is tight by Theorem \ref{P25}.
\item[(iv)] Also, Datta and Singh found \cite{ds} a 21-vertex triangulation of $(\TPSSF)^{\#8}$ and a 26-vertex
triangulation of $(\TPSSF)^{\#14}$ in ${\cal W}_1^{\ast}(4)$. Both are $\mathbb{Z}_2$-tight by Theorem \ref{P25}.
\end{enumerate}
\item[{\bf (e)}] Lutz constructed \cite{lu1} two 12-vertex triangulations of $S^{\,2} \times S^{\,3}$ in ${\cal
W}_2^{\ast}(5)$. Both are tight by Theorem \ref{P25}.

\item[{\bf (f)}] Only finitely many $2k$-dimensional $(k+ 1)$-neighbourly triangulated closed manifolds are known
for $k \geq 2$. By Lemma \ref{L10}, they are all tight. These examples are\,:
 \begin{enumerate}
 \item[(i)] The 9-vertex triangulation of $\mathbb{C}\mathbb{P}^{\,2}$ due to K\"{u}hnel \cite{kb},
 \item[(ii)] six 15-vertex triangulations of homology $\mathbb{H}\mathbb{P}^{\,2}$ (three due to Brehm and
K\"{u}hnel \cite{bk2} and three more due to Lutz \cite{lu2}),
 \item[(iii)] the 16-vertex triangulation of a $K3$-surface due to Casella and K\"{u}hnel \cite{ck}, and
 \item[(iv)] two 13-vertex triangulations of $S^{\,3} \times S^{\,3}$ due to Lutz \cite{lu1}.
\end{enumerate}

\item[{\bf (g)}] Apart from the above list, we know only two tight triangulated manifolds. These are\,:
\begin{enumerate}

\item[(i)] A 15-vertex triangulation of $(\TPSSF) \# (\mathbb{C}\mathbb{P}^{\,2})^{\#5}$ due to Lutz \cite{lu1}.
It is 2-neighbourly, non-orientable, $\mathbb{Z}_2$-tight and in ${\cal W}_2(4)$.

\item[(ii)] A 13-vertex triangulation of $SU(3)/SO(3)$ due to Lutz \cite{lu1}. It is 3-neighbourly, orientable,
$\mathbb{Z}_2$-tight and in ${\cal W}_3(5)$.
\end{enumerate}
The $\mathbb{Z}_2$-tightness of the last two examples does not follow from the results presented here.
\end{enumerate}
Corollary \ref{C11} implies that all the triangulations in Example \ref{E3}\,(a), (b), (c) and (d) are strongly
minimal. By Theorem 4.2 of \cite{ns2}, the triangulations in Example \ref{E3}\,(f) (i) and (ii) are strongly minimal. 
By Theorem 5.8 of \cite{s09}, the 16-vertex triangulation in Example 3.14 
(f) (iii) is strongly minimal. 
By Theorem 4.4 of \cite{ns},  all the triangulations in Example \ref{E3}\,(e) and (f) are minimal. As
far as we know, the minimality of the triangulations in Example \ref{E3}\,(g) is an open problem. }
\end{eg}

We raise the question of how to get more (tight) triangulations meeting the hypothesis of Theorem \ref{P25}. In
particular, we may ask\,:

\begin{qn} \label{Q6}
{\rm Is there a 20-vertex triangulation of $(S^{\,2} \times S^{1})^{\#12}$ or $(\TPSS)^{\#12}$ in ${\cal
W}_1^{\ast}(3)$ or a 20-vertex triangulation of $(S^{\,3}\times S^{\,2})^{\#13}$ in ${\cal W}_2^{\ast}(5)$\,?}
\end{qn}

We do not know for a fact that, for $1<l< (d-1)/2$, the members of ${\cal K}_l(d)$ (or even of ${\cal W}_l(d)$)
actually attain equality in Conjecture \ref{Conj1.6}. Thus, all parts of this conjecture are wide open for $l>1$.
Notice that, as a consequence of Theorem \ref{P21}, the members of ${\cal W}_l^{\ast}(d)$ do attain equality in
Conjecture \ref{Conj1.7} for $1\leq l < (d-1)/2$. However, we do not know if, more generally, the members of
${\cal K}_l^{\ast}(d)$ attain these equalities for $1\leq l < (d-1)/2$, as conjectured. The $l=1$ case is
unrevealing since in this case ${\cal W}_1(d) = {\cal K}_1(d)$ and ${\cal W}_1^{\ast}(d) = {\cal K}_1^{\ast}(d)$.
Moreover, we suspect (but cannot prove) that Theorem \ref{P25} is valid for ${\cal K}^{\ast}_k(d)$ (in place of
the smaller class ${\cal W}^{\ast}_k(d)$). Thus, the most important question raised by this paper is whether (and
to what extent) its results can be extended from $k$-stellated spheres to $k$-stacked spheres. A good place to
begin this investigation is to address the following\,:

\begin{qn} \label{Q7}
{\rm Is Theorem \ref{P19} (on the sigma-vector of $k$-stellated spheres) valid for $k$-stacked spheres\,?}
\end{qn}

\medskip

\noindent {\bf Acknowledgements:} The authors thank the anonymous referee for some useful comments. The second
author was partially supported by grants from UGC Centre for Advanced Study.


\end{document}